\documentclass[a4paper]{article}
\usepackage{amssymb}
\usepackage{amsfonts}
\usepackage{amsmath}
\usepackage{amsthm}
\usepackage{authblk}
\usepackage[USenglish]{babel}
\usepackage{esint}
\usepackage{graphicx}
\usepackage[colorlinks=true]{hyperref}
\usepackage[latin1]{inputenc}

\newtheorem{theorem}{Theorem}[section]
\newtheorem{lemma}[theorem]{Lemma}
\newtheorem{proposition}[theorem]{Proposition}
\newtheorem{corollary}[theorem]{Corollary}
\newtheorem{remark}[theorem]{Remark}
\newtheorem{definition}[theorem]{Definition}
\newtheorem{example}[theorem]{Example}
\newcommand{\N}{\mathbb N}
\newcommand{\Z}{\mathbb Z}

\newcommand{\R}{\mathbb R}

\newcommand{\mbf}{\mathbf}
\newcommand{\mcal}{\mathcal}

\newcommand{\mrm}{\mathrm}

\renewcommand{\a}{\alpha}
\renewcommand{\b}{\beta}
\newcommand{\g}{\gamma}

\renewcommand{\d}{\delta}
\newcommand{\D}{\Delta}

\renewcommand{\t}{\theta}
\newcommand{\Th}{\Theta}
\newcommand{\la}{\lambda}

\newcommand{\s}{\sigma}

\newcommand{\ph}{\varphi}
\renewcommand{\o}{\omega}
\renewcommand{\O}{\Omega}
\newcommand{\wt}{\widetilde}

\newcommand{\ol}{\overline}

\newcommand{\ub}{\underbrace}
\newcommand{\fr}{\frac}
\newcommand{\pa}{\partial}
\newcommand{\n}{\nabla}
\newcommand{\fa}{\forall}

\newcommand{\wk}{\rightharpoonup}
\newcommand{\inc}{\hookrightarrow}

\newcommand{\us}{\underset}

\newcommand{\sm}{\setminus}

\newcommand{\sub}{\subset}
\newcommand{\Sub}{\Subset}

\newcommand{\nin}{\not\in}
\newcommand{\eq}{\equiv}

\newcommand{\x}{\times}

\newcommand{\cd}{\cdot}
\newcommand{\ds}{\dots}

\newcommand{\tx}{\text}
\newcommand{\q}{\quad}
\renewcommand{\l}{\left}
\renewcommand{\r}{\right}

\newcommand{\bthm}{\begin{theorem}}
\newcommand{\ethm}{\end{theorem}}
\newcommand{\blem}{\begin{lemma}}
\newcommand{\elem}{\end{lemma}}
\newcommand{\bprop}{\begin{proposition}}
\newcommand{\eprop}{\end{proposition}}
\newcommand{\bcor}{\begin{corollary}}
\newcommand{\ecor}{\end{corollary}}
\newcommand{\bdefi}{\begin{definition}}
\newcommand{\edefi}{\end{definition}}
\newcommand{\bpf}{\begin{proof}}
\newcommand{\epf}{\end{proof}}
\newcommand{\bl}{\begin{array}{l}}
\newcommand{\bll}{\begin{array}{ll}}
\newcommand{\barr}{\begin{array}}
\newcommand{\earr}{\end{array}}
\newcommand{\bite}{\begin{itemize}}
\newcommand{\eite}{\end{itemize}}
\newcommand{\bequ}{\begin{equation}}
\newcommand{\eequ}{\end{equation}}
\newcommand{\beqa}{\begin{eqnarray}}
\newcommand{\eeqa}{\end{eqnarray}}
\newcommand{\beqy}{\begin{eqnarray*}}
\newcommand{\eeqy}{\end{eqnarray*}}

\allowdisplaybreaks

\begin{document}

\everymath{\displaystyle}

\title{A unified approach of blow-up phenomena for two-dimensional singular Liouville systems}
\author{Luca Battaglia\thanks{Sapienza Universit\`a di Roma, Dipartimento di Matematica, Piazzale Aldo Moro $5$, $00185$ Roma - battaglia@mat.uniroma1.it}, Angela Pistoia\thanks{Sapienza Universit\`a di Roma, Dipartimento di Scienze di Base e Applicate, Via Antonio Scarpa $16$, $00161$ Roma - angela.pistoia@uniroma1.it}}
\date{}

\maketitle\

\begin{abstract}
\noindent We consider generic $2\x2$ singular Liouville systems
\bequ\label{plambda}
\l\{\bll-\D u_1=2\la_1 e^{u_1}-a\la_2 e^{u_2}-2\pi (\a_1-2)\delta_0&\tx{in }\O\\
-\D u_2=2\la_2 e^{u_2}-b\la_1 e^{u_1}-2\pi (\a_2-2)\delta_0&\tx{in }\O\\
u_1=u_2=0&\tx{on }\pa\O\earr\r., 
\eequ
where $\O\ni0$ is a smooth bounded domain in $\mathbb R^2$ possibly having some symmetry with respect to the origin, $\delta_0$ is the Dirac mass at $0,$ $\la_1,\la_2$ are small positive parameters and $a,b,\a_1,\a_2>0$.\\
We construct a family of solutions to \eqref{plambda} which blow up at the origin as $\la_1 \to 0$ and $\la_2 \to 0 $ and   whose  local mass  at the origin  is  a given quantity depending on $a,b,\a_1,\a_2$.\\
In particular, if $ab<4$ we get finitely many possible blow-up values of the local mass, whereas if $ab\ge4$ we get infinitely many.  The blow-up values are produced using an explicit formula which involves Chebyshev polynomials.
\end{abstract}\

%\subjclass[2000]{35J91, 35A01, 35B44, 35B30}
%\keywords{Toda system, blow-up values, tower of bubbles} 

\section{Introduction}\

In this paper we consider the system of singular Liouville equations
\bequ\label{system}
\l\{\bll-\D u_1=2\la_1 e^{u_1}-a\la_2 e^{u_2}-2\pi (\a_1-2)\delta_0&\tx{in }\O\\
-\D u_2=2\la_2 e^{u_2}-b\la_1 e^{u_1}-2\pi (\a_2-2)\delta_0&\tx{in }\O\\
u_1=u_2=0&\tx{on }\pa\O\earr\r., 
\eequ
where $\O\subset\mathbb R^2$ is a smooth bounded domain which contains the origin, $\delta_0$ is the Dirac mass at $0,$ $\lambda_i$ are small positive parameters and the matrix 
\bequ\label{mat-a} \mathcal A=\left(\begin{matrix}2 & -a \\ - b& 2\end{matrix}\right)\ \hbox{ with $a,b>0.$}\eequ
Liouville systems find applications in many fields of physics and mathematics, like theory of chemotaxis \cite{cp}, theory of charged particle beams \cite{kl}, theory of semi-conductors \cite{moc}, Chern-Simons theory \cite{lez,ggo,man,djpt,lezsav,yan1,yan2}, holomorphic projective curves \cite{cal,cw,bjrw,lezsav,bw,dol,gue}.\\

It is not hard to see that any more general Liouville systems
$$\l\{\bll-\D u_1=a_{11}\la_1 e^{u_1}+a_{12}\la_2 e^{u_2}-2\pi (\a_1-2)\delta_0&\tx{in }\O\\
-\D u_2=a_{21}\la_2 e^{u_2}+a_{21}\la_1 e^{u_1}-2\pi (\a_2-2)\delta_0&\tx{in }\O\\
u_1=u_2=0&\tx{on }\pa\O\earr\r.$$
with $a_{12},a_{21}<0<a_{11},a_{22}$ can be brought back to \eqref{system} just by a rescaling of the parameters $\la_1,\la_2$.\\
Using Green's function
\bequ\label{g}
\l\{\bll
-\Delta G(\cdot,y)=\delta_y& \hbox{in}\ \Omega\\
 G(\cdot,y)=0& \hbox{on}\ \partial\Omega
\earr\r.,
\eequ
and its decomposition
\bequ\label{robin}
G(x,y)=\frac 1{2\pi}\log |x-y|+H(x,y)\q\q x,y\in\Omega
\eequ
with $H(x,y)$ smooth, we can eliminate the singularity on the right
hand side of \eqref{system} and rewrite the system as

\bequ\label{system2}
\l\{\bll-\D u_1=2\la_1h_1 e^{u_1}-a\la_2h_2 e^{u_2}&\tx{in }\O\\
-\D u_2=2\la_2h_2 e^{u_2}-b\la_1h_1 e^{u_1}&\tx{in }\O\\
u_1=u_2=0&\tx{on }\pa\O\earr\r., 
\eequ
where 
\bequ\label{accaj}h_i(x)=|x|^{\a_i-2}e^{-2\pi (\a_i-2)H(x,0)}\ \hbox{for}\ i=1,2\eequ
and $e^{-2\pi (\a_i-2)H(x,0)}$ is smooth and positive.

One of the most important and challenging issues concerning Liouville systems \eqref{system} or \eqref{system2} is the blow-up phenomena. A point $x_0\in\ol\O$ is called a blow-up point if a sequence of solutions $u_n=\l(u_{1,n}, u_{2,n }\r)$ satisfies
$$\max\limits_{i=1,2}\max\limits_{B_r(x_0)\cap\O}u_{i,n}=\max\limits_{i=1,2}u_{i,n}(x_n)\us{n\to+\infty}\to+\infty\q \hbox{and}\q x_n\us{n\to+\infty}\to x_0.$$
Knowing the asymptotic behavior of blowing-up solutions near the blow-up points is the first step in applying topological or variational methods to get solutions to the Liouville systems.
In particular, the first main issue is to determine the set of critical masses of solutions with bounded energy, i.e. $\max\limits_{i=1,2}\lambda _{i,n} \int\limits_\Omega h_i e^{u_{i,n}}\le C $ for some $C.$\\
We define the local masses at the blow-up point $x_0$ as
\bequ\label{masses}
m_i(x_0):=\lim\limits_{r\to0}\lim\limits_{n\to+\infty}\lambda_{i,n}\int\limits_{B_r(x_0)}h_ie^{u_{i,n}}\ \hbox{for}\ i=1,2.
\eequ
The local masses have been widely studied in the last years.\\

When the system \eqref{system} reduces to a single singular Liouville equation
\bequ\label{singola}
\l\{\bll-\D u=\la e^u-2\pi (\a-2)\delta_0&\tx{in }\O\\
u=0&\tx{on }\pa\O\earr\r., 
\eequ
 the local mass has been completely characterized. In particular, in the regular case, i.e. $\a=2$, all the blow-up points are internal to $\O$, they are simple and the local mass equals $8\pi$ (see Brezis and Merle \cite{bremer}, Nagasaki and Suzuki \cite{ns}, Li and Shafrir \cite{ls}). In fact, in
this case there is only one bubbling profile: after some rescaling, the bubble approaches a solution
of the Liouville equation
\bequ\label{liou}
\l\{\bll
-\Delta U=e^U& \hbox{in}\ \mathbb R^2,\\
\int\limits_{\mathbb R^2} e^U<+\infty.
\earr\r.
\eequ
 In the singular case, i.e. $\a\not=2$, the local mass around the origin is $4\pi\a $ and the corresponding bubbling profile, after some rescaling, is given by solution to the singular Liouville equation
\bequ\label{liou-s}
\l\{\bll
-\Delta U=|\cd|^{\a-2}e^U& \hbox{in}\ \mathbb R^2,\\
\int\limits_{\mathbb R^2} |\cd|^{\a-2}e^U<+\infty.
\earr\r.
\eequ
 (see Bartolucci and Tarantello \cite{bt} and Bartolucci, Chen, Lin and Tarantello \cite{bclt}).\\
 The knowledge of the bubbling profile is the main step in finding existence and multiplicity results concerning the equation \eqref{singola}. Indeed, bubbling solutions with multiple concentration points have been built by Baraket and Pacard \cite{barpac}, del Pino, Kowalczyk and Musso \cite{dkm} and Esposito,
Grossi and Pistoia \cite{egp} in the regular case and by del Pino, Esposito and Musso \cite{dem} and D'Aprile \cite{d} in the singular case. Moreover, a degree formula has been obtained by Chen and Lin \cite{cl1,cl2} and Malchiodi \cite{mal}, whereas solutions have also been found through variational methods by Bartolucci and Malchiodi
\cite{barmal}, Bartolucci, De Marchis and Malchiodi \cite{bdm}, Carlotto and Malchiodi \cite{carmal}, Djadli \cite{dja} and Malchiodi and Ruiz \cite{mr1}. 
 \\
 
 The natural generalization to \eqref{singola} is the $2 \x 2$ system \eqref{system} when the matrix $\mathcal A=(a_{ij})_{2\x2}$ is as in \eqref{mat-a}. In particular, when $\mathcal A$ is the Cartan matrix of a simple Lie algebra we get the well-known Toda system. In this case,
 since the rank of the simple Lie Algebra is $2$, there are three types of corresponding
Cartan matrices of rank $2$:
 \bequ\label{lie}
 \mathcal A_2=\left(\begin{matrix}2 & -1 \\ -1 & 2\end{matrix}\right),\quad \mathcal B_2=\left(\begin{matrix}2 & -1 \\ -2 & 2\end{matrix}\right), \quad \mathcal G_2=\left(\begin{matrix}2 & -1 \\ -3 & 2\end{matrix}\right). 
 \eequ
 In the regular $\mathcal A_2-$Toda system, i.e. $\a_1=\a_2=2,$ Jost, Lin and Wang in \cite{jlw} found that the local masses can only take $5$ values. Moreover, all these blow-up values can occur as shown by Musso, Pistoia and Wei in \cite{mpw} (see also Ao and Wang in \cite{aw}).\\
In the singular case Lin, Wei and Zhang in \cite{lwzhang} found that only $5$ possible values are allowed for the local masses, provided the singularities $\a_1$ and $\a_2$ satisfy a suitable condition (which also include the regular case) (see Example \ref{a2}).\\
Recently, Lin and Zhang in \cite{lz} found that only $7$ possible values are allowed for the local masses in the regular $\mathcal B_2-$Toda system (see Example \ref{b2}) and $11$ possible values are allowed for the local masses in the regular $\mathcal G_2-$Toda system (see Example \ref{g2}) under some extra assumptions.\\
Solutions to the regular $A_2$ Toda system have been found both through the computation of the degree by Lin, Wei and Yang in  \cite{lwyang} and variationally by Battaglia, Jevnikar, Malchiodi and Ruiz in \cite{bjmr}, Jevnikar, Kallel and Malchiodi \cite{jkm}, Malchiodi and Ndiaye \cite{mn} and Malchiodi and Ruiz \cite{mr2}. Variational solutions have also been found for the $A_2$ Toda system in Battaglia, Jevnikar, Malchiodi and Ruiz in \cite{bjmr}, Battaglia in \cite{bat} and Battaglia and Malchiodi in  \cite{batmal} and for the $B_2$ and $G_2$ systems by Battaglia in  \cite{bat2}.\\

At this stage, two questions naturally ariese:
\bite
\item[(Q1)] {\em which are the values of the local masses at the origin for the system \eqref{system} for a general matrix $\mathcal A$ with or without singular sources?}
\item[(Q2)] {\em are these values attained?}
\eite

In this paper we focus on the second question and we give a partial answer. More precisely, we build solutions to the system \eqref{system}, whose components blows-up at the origin and whose local masses are quantized in
terms of $a,$ $b$, $\a_1$ and $\a_2.$ In particular, if $\det(\mathcal A)\le 0$ we find infinitely many possible values for the local masses.
We also provide an explicit formula involving Chebyshev polynomials which produces  blow-up values of the local masses (see Remark \ref{ricetta}).  We also conjecture that these   are the only admissible values when the blowing-up profile of each component resembles one or more bubbles solving the scalar Liouville equations \eqref{singola}. Indeed, they coincide  with the known ones when the matrix $\mathcal A$ is as in \eqref{lie} (see Examples \ref{a2}, \ref{b2} and \ref{g2}).
  \\

Let us state our main result.\

For any integer $\ell\in\mathbb N$ we introduce the polynomials 
\bequ\label{pj}
\l\{\bll
P_0(t)=0\\
P_1(t)=1\\
P_2(t)=1\\
P_\ell(t)=\prod_{i=1}^{\l[\fr{\ell-1}2\r]}\l(t-2-2\cos\fr{2\pi i}\ell\r)& \hbox{if}\ \ell\ge3
\earr\r..
\eequ
and the real numbers $\beta_\ell=\beta_\ell(a,b,\a_1,\a_2)$ defined as follows
\bequ\label{betaj}
\b_{\ell} = \left\{\bll\a_1P_{\ell}(ab)+a\a_2P_{\ell-1}(ab)&\hbox{if $\ell$ is odd}\\
b\a_1P_{\ell}(ab)+\a_2P_{\ell-1}(ab) & \hbox{if $\ell$ is even}.\earr\right.
\eequ
 Then we define the (possibly infinite) integer
\bequ\label{kmax}
k_{\max}=k_{\max}(a,b,\a_1,\a_2):=\sup\{k:\,\b_\ell>0,\,\fa\,\ell=1,\ds,k\}.
\eequ
By \eqref{pj} and \eqref{betaj} it immediately follows that $k_{\max}\ge2$, and we also deduce that 
$k_{\max}=+\infty$ if $ab\ge4.$ In Remark \ref{tmax3}, we find the following expression of $k_{\max}$ in terms of $a$ and $b$ when $ab<4$
\bequ\label{kmassimo}
k_{\max}=\l\{\bll\fr{2\pi}{\arccos\l(\fr{ab}2-1\r)}&\tx{if }\fr{2\pi}{\arccos\l(\fr{ab}2-1\r)}\in\N\\\l[\fr{2\pi}{\arccos\l(\fr{ab}2-1\r)}\r]&\tx{if }\b_{\l[\fr{2\pi}{\arccos\l(\fr{ab}2-1\r)}\r]+1}<0\\\l[\fr{2\pi}{\arccos\l(\fr{ab}2-1\r)}\r]+1&\tx{if }\b_{\l[\fr{2\pi}{\arccos\l(\fr{ab}2-1\r)}\r]+1}>0 \earr\r..\eequ

\begin{definition}\label{omega} ${}$\\
Set $\mathcal I:=\{\ell\in\{1,\ds,k\}\ :\ {\b_\ell} \in 2\N\}.$  
We say that $\O$ is  \emph{compatible}   if  
$$e^{\fr{2\pi}{\mathfrak m }\iota}\O:=\l\{\l(x_1\cos \fr{2\pi}{\mathfrak m } -x_2\sin \fr{2\pi}{\mathfrak m } ,x_1\sin \fr{2\pi}{\mathfrak m } +x_2\cos \fr{2\pi}{\mathfrak m }\r)\  :\ (x_1,x_2)\in\O\r\}=\O ,$$
 where ${\mathfrak m }:=l.c.m.\left\{{{\mathfrak m }_\ell}\in\mathbb N\ :\ {\b_\ell\over  {\mathfrak m }_\ell}\not\in 2\mathbb N,\ \ell\in\mathcal I\right\}.$\ 
In particular, if $\mathcal I=\emptyset$ any smooth bounded domain $\Omega$ containing the origin is  \emph{compatible}.
%Let $\b_1,\ds,\b_k$ be given. We say that $\O$ is \emph{compatible} if there exists $m\in\N$ such that $\fr{\b_\ell}{2m}\nin 2\N$ for any $\ell=1,\ds,k$ and
%$$e^{\fr{2\pi}m\iota}\O:=\l\{\l(\cos\l(\fr{2\pi}m\r)x_1-\sin\l(\fr{2\pi}m\r)x_2,\sin\l(\fr{2\pi}m\r)x_1+\cos\l(\fr{2\pi}m\r)x_2\r):\,(x_1,x_2)\in\O\r\}=\O$$
% particular, if $\b_\ell\nin2\N$ for any $\ell$, then any smooth bounded domain $\Omega$ %containing the origin is \emph{compatible}, since the previous condition is trivially %satisfied by $m=1$.
\end{definition}\

\begin{remark}${}$\\
The integer $\mathfrak m$ introduce in the previous definition is not uniquely defined, since it depends on the choice of $\mathfrak m_1,\ds,\mathfrak m_k$, which are not unique. In the definition of compatibility we want the equality to hold true for \emph{at least one} of such possible $\mathfrak m$'s.
\end{remark}\

\bthm\label{main}${}$\\
Let $k\in\N,\,k\le k_{\max}$ be fixed, $\b_1,\ds,\b_k$ be defined by \eqref{betaj} and $\O\ni0$ a smooth bounded domain which is \emph{compatible} in the sense of Definition \ref{omega}.\\
Then, there exists $\ol\la=\ol\la(k)>0$ such that for $\la$ satisfying
\bequ\label{cond}\l\{\bll\la_1,\la_2\in\l(0,\ol\la\r)&\tx{if }k<k_{\max}\\
\la_1,\la_2\in\l(0,\ol\la\r),\,\la_2\le\la_1^\fr{\g-\b_{k_{\max}+1}}{\b_{k_{\max}}}&\tx{if }k=k_{\max}\tx{ is odd}\\
\la_1,\la_2\in\l(0,\ol\la\r),\,\la_1\le\la_2^\fr{\g-\b_{k_{\max}+1}}{\b_{k_{\max}}}&\tx{if }k=k_{\max}\tx{ is even}\earr\r.\q\tx{for some }\g>0.\eequ
the problem \eqref{plambda} has a solution $u=u_\la=(u_{1,\la},u_{2,\la})$.\\
Moreover, there holds
\bequ
\label{mass}
 m_1(0)=2\pi\sum_{j=0}^{\l[\fr{k-1}2\r]}\b_{2j+1}\quad\quad \hbox{and}\quad\quad
 m_2(0) =2\pi\sum_{j=0}^{\l[\fr{k-2}2\r]}\b_{2j+2},\eequ
 where we agree that $m_2(0)=0$ if $k=1.$
Moreover, if $G$ is the Green's function defined in \eqref{g}, we have as $\lambda\to0$
\bequ
\label{green}u_1 \to[2m_1(0)-am_2(0)]G(\cd,0)\quad\quad \hbox{and}\quad\quad u_2 \to[2m_2(0)-bm_1(0)]G(\cd,0)
\eequ
weakly in $W^{1,q}(\O)$ for any $q<2$ and strongly in $C^\infty_{\mrm{loc}}(\O\sm\{0\})$.
\ethm\

\begin{remark}${}$\\
We point out that if $\Omega$ is a symmetric domain according to in Definition \ref{omega}, then the functions $h_j$ defined in \eqref{accaj} and  the solutions found in Theorem \ref{main} inherit the symmetry of the domain $\O$, namely they satisfy the symmetry condition
  $u\l(e^{\fr{2\pi}{\mathfrak m }\iota}x\r)=u(x)$ for any $x\in\Omega$, where ${\mathfrak m }$ is as in Definition \ref{omega}.
 \end{remark}

\begin{remark}\label{ricetta}${}$\\
As far as we know, this is the first result which gives a clear relation between the local masses at the origin and their possible number and the values of the entries of the matrix $\mathcal A$ in \eqref{mat-a} and the values of the singularities $\a_1$ and $\a_2.$ Indeed, 
we can express the masses in \eqref{mass} in terms of the value of the polynomials $P_\ell(t)$ at $t=ab$ as (see Remark \ref{masse})
\bequ
\label{mass1} m_1(0)= \l\{\bll2\pi aP_{\l[\fr{k-1}2\r]}(ab)\l(b\a_1P_{\l[\fr{k-1}2\r]}(ab)+\a_2P_{\l[\fr{k-3}2\r]}(ab)\r)&\tx{if }k\in(4\N+1)\cup(4\N+2)\\2\pi P_{\l[\fr{k-1}2\r]}(ab)\l(\a_1P_{\l[\fr{k-1}2\r]}(ab)+a\a_2P_{\l[\fr{k-3}2\r]}(ab)\r)&\tx{if }k\in(4\N+3)\cup4\N\earr\r.\eequ
and
\bequ
\label{mass2} m_2(0)= \l\{\bll2\pi P_{\l[\fr{k-2}2\r]}(ab)\l(b\a_1P_{\l[\fr{k}2\r]}(ab)+\a_2P_{\l[\fr{k-2}2\r]}(ab)\r)&\tx{if }k\in4\N\cup(4\N+1)\\2\pi bP_{\l[\fr{k-2}2\r]}(ab)\l(\a_1P_{\l[\fr{k}2\r]}(ab)+a\a_2P_{\l[\fr{k-2}2\r]}(ab)\r)&\tx{if }k\in(4\N+2)\cup(4\N+3)\earr\r.
\eequ
where the range of $k$ is between $1$ and the number $k_{\max}$ defined in \eqref{kmax}.
\end{remark}

\begin{remark}${}$\\
The bubbling profile of each component resembles a sum (with alternating sign) of bubbles solutions to different singular Liouville problems \eqref{liou-s}: all the bubbles are centered at the origin and the rate of concentration of each bubble at the origin is slower than the previous one, namely
 \bequ\label{raf-ans}u_1\sim w_1-\fr{a}2 w_2+ w_3-\frac a2 w_4+\dots\qquad \hbox{and}\qquad u_2\sim -\fr{b}2 w_1+ w_2-\frac b2 w_3+ w_4+\dots\eequ
where 
\bequ\label{wi}w_i(x):=\log2\b_i^2\fr{\d_i^{\b_i}}{\l(\d_i^{\b_i}+|x|^{\b_i}\r)^2}\quad x\in\mathbb R^2,\ \delta_i>0\quad \hbox{solves}\quad
-\Delta w_i=|\cd|^{\b_i-2}e^{w_i}\ \hbox{in}\ \mathbb R^2\eequ
and $\fr{\d_i}{\d_{i+1}}$ approaches zero.
The construction of a solution with such a profile is possible as long as the exponents $\b_i$'s are positive and that is why we need to introduce the maximal number of bubbles $k_{\max}$ in \eqref{kmax}. Moreover, each bubble $w_\ell$ scaled with $\d_i$ turns out to be a singular source for the equation solved by the bubble $w_i$ whenever $\ell<i$.\\
Therefore, the choice of each $\b_i$ takes into account the singular sources present in the equation and all the singular sources generated by the interactions between the bubbles $w_i$ and all the previous ones. This fact leads to choose $\b_i$ as in \eqref{beta-rec} to ensure that the prescribed profile is {\em almost} a solution to system \eqref{system} (as carefully proved in Lemma \ref{theta} and Lemma \ref{rla}).\\
 This kind of construction is strongly inspired by the bubble-tower construction
in Musso, Pistoia and Wei \cite{mpw} (see also Grossi and Pistoia \cite{gp}), where the regular $\mathcal A_2-$Toda system was studied. Nevertheless, the general case turns out to be rather
delicate.\\
In particular, the interaction between the two components is much more involved because the concentration of each bubble is affected by all the other previous bubbles, not only the ones for which the same component concentrates. Even and odd bubbles affect the concentration in opposite ways.\\
Moreover, we will need some rather involved symmetry condition, which are needed to invert a linearized operator and strongly depend on the values of $\b_i$. Finally, the presence of singularites gives weaker regularity properties and makes some estimates more subtle.
\end{remark}

 In the following examples we describe how our result can be applied to classical problems.

\begin{example}\label{a2}${}$\\
If $a=b=1$, the system \eqref{system2} becomes the well-known \underline{$A_2$-Toda system}
$$\l\{\bll-\D u_1=2\la_1h_1e^{u_1}-\la_2h_2e^{u_2}&\tx{in }\O\\
-\D u_2=2\la_2h_2e^{u_2}-\la_1h_1e^{u_1}&\tx{in }\O\\
u_1=u_2=0&\tx{on }\pa\O\earr\r.$$
We have $ab=1$ and by \eqref{kmassimo} we compute $k_{\max}=3.$ Moreover, 
$$\l\{\bl
P_1(1)=1\\
P_2(1)=1\\
P_3(1)=-1-2\cos{2\pi\over3}=0
\earr\r..$$
Then, by \eqref{betaj} and \eqref{mass} (possibly exchanging the role of the components) we deduce
 the following configurations for $(m_1(0),m_2(0)):$
 \bite
 \item if $k=1$ we get $2\pi(\alpha_1,0)$ and $2\pi(0, \alpha_2),$ 
 \item if $k=2$ we get $2\pi(\a_1,\a_1+\a_2)$ and $2\pi(\a_1+\a_2,\a_2),$
 \item if $k=3$ we get $2\pi(\a_1+\a_2,\a_1+\a_2) .$
 \eite
 
In \cite{lwzhang} Lin, Wei and Zhang show that, for suitable values of $\a_1,\a_2$ (including the regular case $\a_1=\a_2=2$), the only possible values are the five above. Therefore, Theorem \ref{main} shows in particular the sharpness of their classification.\\
For the regular Toda system, Theorem \ref{main} was already proved by Musso, Pistoia and Wei \cite{mpw}.
\end{example}\

\begin{example}\label{b2}${}$\\
The case $a=1,\,\a_1=\a_2=b=2$ is the \underline{$B_2$-Toda system}
$$\l\{\bll-\D u_1=2\la_1e^{u_1}-\la_2e^{u_2}&\tx{in }\O\\
-\D u_2=2\la_2e^{u_2}-2\la_1e^{u_1}&\tx{in }\O\\
u_1=u_2=0&\tx{on }\pa\O\earr\r. $$

We have $ab=2$ and by \eqref{kmassimo} we compute $k_{\max}=4.$ Moreover, 
$$\l\{\bl
P_1(2)=1\\
P_2(2)=1\\
P_3(2)=-2\cos{2\pi\over3}=1\\
P_4(2)=-2\cos{ \pi\over2}=0
\earr\r..$$
Then, by \eqref{betaj} and \eqref{mass} we deduce
 the following configurations for $(m_1(0),m_2(0)):$
 \bite
 \item if $k=1$ we get $2\pi (\beta_1,0)=2\pi(2,0)$ 
 \item if $k=2$ we get $2\pi ( \beta_1, \beta_2)=2\pi(2,6)$
 \item if $k=3$ we get $2\pi(\beta_1+\beta_3,\beta_2)=2\pi(6,6)$
 \item if $k=4$ we get $2\pi(\beta_1+\beta_3,\beta_2+\beta_4)=2\pi(6,8)$
 \eite
 and exchanging the role of the components (i.e. $b=1$ and $a=2$)
 \bite
 \item if $k=1$ we get $2\pi (0,\beta_1)=2\pi(0,2)$ 
 \item if $k=2$ we get $2\pi ( \beta_2, \beta_1)=2\pi(4,2)$
 \item if $k=3$ we get $2\pi(\beta_2,\beta_1+\beta_3)=2\pi(4,8)$
 \item if $k=4$ we get $2\pi(\beta_2+\beta_4,\beta_1+\beta_3)=2\pi(6,8)$
 \eite 
In \cite{lz} Lin and Zhang show that no other values are admissible in case of blow up. Theorem \ref{main} shows the sharpness of their classification.\\\\
\end{example}\

\begin{example}\label{g2}${}$\\
The case $a=1,\,b=3,\,\a_1=\a_2=2$ is the \underline{$G_2$-Toda system}
$$ \l\{\bll-\D u_1=2\la_1e^{u_1}-\la_2e^{u_2}&\tx{in }\O\\
-\D u_2=2\la_2e^{u_2}-3\la_1e^{u_1}&\tx{in }\O\\
u_1=u_2=0&\tx{on }\pa\O\earr\r..$$
We have $ab=3$ and by \eqref{kmassimo} we compute $k_{\max}=6.$ Moreover, 
$$\l\{\bl
P_1(3)=1\\
P_2(3)=1\\
P_3(3)=1-2\cos{2\pi\over3}=2\\
P_4(3)=1-2\cos{ \pi\over2}=1\\
P_5(3)=\left(1-2\cos{2\pi\over5}\right)\left(1-2\cos{4\pi\over5}\right)=1\\ P_6(3)=\left(1-2\cos{\pi\over3}\right)\left(1-2\cos{2\pi\over3}\right)=0
\earr\r..$$
Then, by \eqref{betaj} and \eqref{mass} we deduce
 the following configurations for $(m_1(0),m_2(0)):$
 \bite
 \item if $k=1$ we get $2\pi (\beta_1,0)=2\pi(2,0)$,
 \item if $k=2$ we get $2\pi ( \beta_1, \beta_2)=2\pi(2,8)$,
 \item if $k=3$ we get $2\pi(\beta_1+\beta_3,\beta_2)=2\pi(8,8)$,
 \item if $k=4$ we get $2\pi(\beta_1+\beta_3,\beta_2+\beta_4)=2\pi(8,18)$,
 \item if $k=5$ we get $2\pi(\beta_1+\beta_3+\beta_5,\beta_2+\beta_4)=2\pi(12,18)$,
 \item if $k=6$ we get $2\pi(\beta_1+\beta_3+\beta_5,\beta_2+\beta_4+\beta_6)=2\pi(12,20)$,
 \eite
 and exchanging the role of the components (i.e. $b=1$ and $a=3$);
 \bite
 \item if $k=1$ we get $2\pi (0,\beta_1)=2\pi(0,2)$,
 \item if $k=2$ we get $2\pi ( \beta_2, \beta_1)=2\pi(4,2)$,
 \item if $k=3$ we get $2\pi(\beta_2,\beta_1+\beta_3)=2\pi(4,12)$,
 \item if $k=4$ we get $2\pi(\beta_2+\beta_4,\beta_1+\beta_3)=2\pi(10,12)$,
 \item if $k=5$ we get $2\pi(\beta_2+\beta_4,\beta_1+\beta_3+\beta_5)=2\pi(10,20)$,
 \item if $k=6$ we get $2\pi(\beta_2+\beta_4+\beta_6,\beta_1+\beta_3+\beta_5)=2\pi(12,20)$.
 \eite 
In \cite{lz} Lin and Zhang found the previous blow-up values under some extra assumptions. Theorem \ref{main} shows that these blow-up values are attained. \\
\end{example}\

\begin{example}\label{seno}${}$\\
The case $ab=4$ is particularly interesting: in Remark \ref{tmax2} this is the borderline scenario to have an infinite number of blow-up values. This fact is related to the matrix of the coefficients in \eqref{plambda} being singular.\\
In fact, if we consider the system
$$\l\{\bll-\D u_1=2\la_1h_1e^{u_1}-a\la_2h_2e^{u_2}&\tx{in }\O\\
-\D u_2=2\la_2h_2e^{u_2}-\fr{4}a\la_1h_1e^{u_1}&\tx{in }\O\\
u_1=u_2=0&\tx{on }\pa\O\earr\r.,$$
then a suitable linear combination of the two equation gives $\l\{\bll-\D\l(u_1+\fr{a}2u_2\r)=0&\tx{in }\O\\
u_1+\fr{a}2u_2=0&\tx{on }\pa\O\earr\r.$, which means $u_2=-\fr{2}au_1$; therefore, in this case \eqref{system2} is equivalent to the scalar equation
$$\l\{\bll-\D u=2\la_1h_1e^{u}-a\la_2h_2e^{-\fr{2}au}&\tx{in }\O\\
u=0&\tx{on }\pa\O\earr\r..$$
In this case (see Remark \ref{t=4})
$$P_\ell(4)=\prod_{i=1}^{\l[\fr{\ell-1}2\r]}2\l(1- \cos\fr{2\pi i}\ell\r)=\left\{\begin{aligned}
&\ell\ \hbox{if $\ell$ is odd}\\ &\frac \ell 2\ \hbox{if $\ell$ is even}. \end{aligned}\right.$$

Therefore, using Remark \ref{masse}, the infinitely many blow-up masses are
$$2\pi(\a_1,0),\q 2\pi\l(\a_1,\fr{4}a\a_1+\a_2\r),\ 2\pi\l(4\a_1+a\a_2,\fr{4}a\a_1+\a_2\r),\ 2\pi\l(4\a_1+a\a_2,\fr{12}a\a_1+4\a_2\r),\ \cdots$$
$$2\pi\l((\ell+1)^2\a_1+\fr{a}2\ell(\ell+1)\a_2,\fr{2}a\ell(\ell+1)\a_1+\ell^2\a_2\r),\ 2\pi\l((\ell+1)^2\a_1+\fr{a}2\ell(\ell+1)\a_2,\fr{2}a(\ell+1)(\ell+2)\a_1+(\ell+1)^2\a_2\r),\ \cdots$$
$$2\pi(0,\a_2),\q 2\pi(\a_1+a\a_2,\a_2),\ 2\pi\l(\a_1+a\a_2,\fr{4}a\a_1+4\a_2\r),\ 2\pi\l(4\a_1+3a\a_2,\fr{4}a\a_1+4\a_2\r),\ \cdots$$
$$2\pi\l(\ell^2\a_1+\fr{a}2\ell(\ell+1)\a_2,\fr{2}a\ell(\ell+1)\a_1+(\ell+1)^2\a_2\r),\ 2\pi\l((\ell+1)^2\a_1+\fr{a}2(\ell+1)(\ell+2)\a_2,\fr{2}a\ell(\ell+1)\a_1+(\ell+1)^2\a_2\r),\ \cdots$$
The case $\a_1=\a_2=a=2$ is known as \underline{Sinh-Gordon equation}. The above-mentioned values are shown to be the only admissible ones for any blow-up, as showed by Jost, Wang, Ye and Zhou in \cite{jwyz}. Moreover, all such values had already proved to be attained by Grossi and Pistoia in \cite{gp}, where Theorem \ref{main} is proved in this particular case.\\
The case $\a_1=\a_2,\,a=1$ is known as \underline{Tzitzeica equation}. Jevnikar and Yang in \cite{jy} proved that no other value, besides the ones above, can occur for blow-up masses.\\
Moreover, for $\a_1=\a_2=2$, the above-mentioned blow-up values are attained for any $a>0$, as Pistoia and Ricciardi have recently showed in \cite{pr}.
\end{example}\

The proof of our result relies on a contraction mapping argument and it is performed in Section \ref{proof}. In Section \ref{theansatz} we give a more precise description of the leading term \eqref{raf-ans}, in Section \ref{errors} we estimate the error terms and in Section \ref{linear} we study the linear theory.\\
The symmetry introduced in the definition \ref{omega} is a technical condition used in the linear theory which ensures the non-degeneracy in a one-codimensional space of the bubble $w_i$ defined by \eqref{wi} even when the parameter $\beta_i$ is even (see the Appendix \ref{appe}).

\section{The ansatz}\label{theansatz}\

For any $\b>0$, let 
$$w^\b_\d(x):=\log2\b^2\fr{\d^\b}{(\d^\b+|x|^\b)^2}\quad x\in\mathbb R^2,\ \delta>0$$
be the solutions to the singular Liouville problem in the whole plane, namely
$$\l\{\bll-\D w^\b_\d=|\cd|^{\b-2}e^{w^\b_\d}&\tx{in }\R^2\\\int_{\R^2}|\cd|^{\b-2}e^{w^\b_\d}<+\infty\earr\r..$$
For any integer $k\in[1,k_{\max}]$ we will look for a solution to problem \eqref{plambda} as 
\bequ\label{ans}u_\la=W_\la+\phi_\la=(W_{1,\la}+\phi_{1,\la},W_{2,\la}+\phi_{2,\la}).\eequ
The components of the main term $W_\la$ are defined as
\bequ\label{mainterm}
\left\{\begin{aligned}
W_{1,\la}&=\mrm Pw_1-\fr{a}2\mrm Pw_2\ds=\sum_{j=0}^{\l[\fr{k-1}2\r]}\mrm Pw_{2j+1}-\fr{a}2\sum_{j=0}^{\l[\fr{k-2}2\r]}\mrm Pw_{2j+2}\\
\nonumber W_{2,\la}&=-\fr{b}2\mrm Pw_1+\mrm Pw_2\ds=\sum_{j=0}^{\l[\fr{k-2}2\r]}\mrm Pw_{2j+2}-\fr{b}2\sum_{j=0}^{\l[\fr{k-1}2\r]}\mrm Pw_{2j+1}\q\q\q 
\end{aligned}\right.\eequ
where we agree that if $k=1$ the second sum in $W_{1,\la}$ and the first sum in $W_{2,\la}$ are zero.

Moreover $w_\ell:=w^{\b_\ell}_{\d_\ell}$ and the projection $\mrm P:H^1(\O)\to H^1_0(\O)$ is defined by
\bequ\label{pro}
 \l\{\bll-\D(\mrm Pu)=-\D u&\tx{in }\O\\\mrm Pu=0&\tx{on }\pa\O\earr\r..
\eequ
The $\b_\ell$'s are defined by recurrence as
\bequ\label{beta-rec}
\left\{\begin{aligned}
\b_1&=\alpha_1\\ \b_2&=b\alpha_1+\alpha_2\\
\b_{2j+1}&=a\sum_{i=0}^{j-1}\b_{2i+2}-2\sum_{i=0}^{j-1}\b_{2i+1}+\a_1=a\b_{2j}-\b_{2j-1}\\
\b_{2j+2}&=b\sum_{i=0}^j\b_{2i+1}-2\sum_{i=0}^{j-1}\b_{2i+2}+\a_2=b\b_{2j+1}-\b_{2j}
\end{aligned}\right..
\eequ 
Actually, the two definitions of $\beta_\ell$'s given in \eqref{betaj} and in \eqref{beta-rec} match perfectly. That will be proved in Section \ref{polinomi}.

The concentration parameters $\delta_\ell$'s satisfy
\beqa
\nonumber-\b_{2j+1}\log\d_{2j+1}-2\sum_{i=j+1}^{\l[\fr{k-1}2\r]}\b_{2i+1}\log\d_{2i+1}+a\sum_{i=j}^{\l[\fr{k-2}2\r]}\b_{2i+2}\log\d_{2i+2}-\log\l(2\b_{2j+1}^2\r)&&\\
\label{deltaj}+2\pi\l(2\sum_{i=0}^{\l[\fr{k-1}2\r]}\b_{2i+1}-a\sum_{i=0}^{\l[\fr{k-2}2\r]}\b_{2i+2}-\a_1+2\r)H(0,0)+\log(2\la_1)&=&0,\\
\nonumber-\b_{2j+2}\log\d_{2j+2}-2\sum_{i=j+1}^{\l[\fr{k-2}2\r]}\b_{2i+2}\log\d_{2i+2}+b\sum_{i=j+1}^{\l[\fr{k-1}2\r]}\b_{2i+1}\log\d_{2i+1}-\log\l(2\b_{2j+2}^2\r)&&\\
\nonumber+2\pi\l(2\sum_{i=0}^{\l[\fr{k-2}2\r]}\b_{2i+2}-b\sum_{i=0}^{\l[\fr{k-1}2\r]}\b_{2i+1}-\a_2+2\r)H(0,0)+\log(2\la_2)&=&0.
\eeqa
The choice of $\beta_\ell$'s and $\delta_\ell$'s is motivated by Lemma \ref{theta}. 

It is useful to point out that by \eqref{deltaj} we easily deduce that
$$\d_{2j+1}=d_{2j+1}\l\{\bll\la_1^\fr{P_{k-2j}(ab)}{\b_{2j+1}}\la_2^\fr{aP_{k-2j-1}(ab)}{\b_{2j+1}}&\tx{if }k\tx{ is odd}\\\la_1^\fr{P_{k-2j-1}(ab)}{\b_{2j+1}}\la_2^\fr{aP_{k-2j}(ab)}{\b_{2j+1}}&\tx{if }k\tx{ is even}\earr\r.\q\hbox{and}\q\d_{2j+2}=d_{2j+2}\l\{\bll\la_1^\fr{bP_{k-2j-1}(ab)}{\b_{2j+2}}\la_2^\fr{P_{k-2j-2}(ab)}{\b_{2j+2}}&\tx{if }k\tx{ is odd}\\\la_1^\fr{bP_{k-2j-2}(ab)}{\b_{2j+2}}\la_2^\fr{P_{k-2j-1}(ab)}{\b_{2j+2}}&\tx{if }k\tx{ is even}\earr\r.,$$
which implies 
$$\fr{\d_\ell}{\d_{\ell+1}}=\fr{d_\ell}{d_{\ell+1}}\l\{\bll\la_1^\fr{\b_{k+1}}{\b_\ell\b_{\ell+1}}\la_2^\fr{\b_k}{\b_\ell\b_{\ell+1}}&\tx{if }k\tx{ is odd}\\\la_1^\fr{\b_k}{\b_\ell\b_{\ell+1}}\la_2^\fr{\b_{k+1}}{\b_\ell\b_{\ell+1}}&\tx{if }k\tx{ is even}\earr\r..$$
We want to have $\fr{\d_\ell}{\d_{\ell+1}}\us{\la\to0}\to0$ for any $\ell$, i.e. each bubble is slower than the previous one; this is always satisfied if $\beta_k+1>0$, namely $k<k_{\max}$, otherwise we need the additional condition in \eqref{cond}. The condition $\eqref{cond}$ also ensures that $\fr{\d_\ell}{\d_{\ell+1}}=O(|\la|^\g)$ for some $\g>0$, which will be useful in some estimates throughout the paper.\\

Finally, the remainder term $\phi_\la$ in \eqref{ans} belongs to the following space
$$\mbf H:=\l\{\phi=(\phi_1,\phi_2)\in H^1_0(\O)\x H^1_0(\O):\,\phi_i\l(e^{\fr{2\pi}{\mathfrak m }\iota}x\r)=\phi_i(x)\ \hbox{for any $x\in\Omega$,\, $i=1,2$}\r\},$$
where ${\mathfrak m }$ is as in Definition \ref{omega}. We agree that if ${\mathfrak m }=1$ than $\mbf H$ is nothing but the space $H^1_0(\O)\x H^1_0(\O)$.\\

The space $H^1_0(\O)\x H^1_0(\O)$ is equipped with the norm
$$ \|(u_1,u_2)\|:=\|u_1\|+\|u_2\|,\q\q\q \hbox{where}\q \|u\|:=\l(\int_\O|\n u|^2\r)^\fr{1}2.$$
Moreover, we also consider the space $L^p(\O)\x L^p(\O),$ with $p>1,$ equipped with the norm
$$\|(u_1,u_2)\|_p:=\|u_1\|_p+\|u_2\|_p,\q\q\q \hbox{where}\q \|u\|_p:=\l(\int_\O|u|^p\r)^\fr{1}p.$$

\subsection{The choice of concentration parameters}
For any integer $\ell=1,\dots,k$ we introduce the function $\Theta_\ell$ which reads if $\ell$ is odd, i.e. $\ell=2j+1$ as

\beqa
\nonumber\Theta_{2j+1}(y)&=&\l(\sum_{i=0}^{\l[\fr{k-1}2\r]}\mrm Pw_{2i+1}-w_{2j+1}-\fr{a}2\sum_{i=0}^{\l[\fr{k-2}2\r]}\mrm Pw_{2i+2}\r)(\d_{2j+1}y)-(\b_{2j+1}-\a_1)\log|\d_{2j+1}y|\\
\label{thetaj}&-&2\pi(\a_1-2)H(\d_{2j+1}y,0)+\log(2\la_1),
\eeqa
and if $\ell$ is even, i.e. $\ell=2j+2$, as
\beqa
\Theta_{2j+2}(y)&=&\l(\sum_{i=0}^{\l[\fr{k-2}2\r]}\mrm Pw_{2i+2}-w_{2j+2}-\fr{b}2\sum_{i=0}^{\l[\fr{k-1}2\r]}\mrm Pw_{2i+1}\r)(\d_{2j+2}y)-(\b_{2j+2}-\a_2)\log|\d_{2j+2}y|\\
\label{thetaj+}&-&2\pi(\a_2-2)H(\d_{2j+2}y,0)+\log(2\la_2).
\eeqa
We agree that if $k=1$ the second sum in \eqref{thetaj} and the first sum in \eqref{thetaj+} are zero.
We shall estimate each functions $\Theta_\ell$ on the corresponding scaled annulus
$$\fr{\mcal A_\ell}{\d_\ell}=\l\{y\in\fr{\O}{\d_\ell}:\,\fr{\sqrt{\d_{\ell-1}\d_\ell}}{\d_\ell}\le|y|\le\fr{\sqrt{\d_\ell\d_{\ell+1}}}{\d_\ell}\r\}, \hbox{where}\ 
\mcal A_\ell:=\l\{x\in\O:\,\sqrt{\d_{\ell-1}\d_\ell}\le|x|\le\sqrt{\d_\ell\d_{\ell+1}}\r\},$$
where we agree that $\d_0=0$ and $\d_{k+1}=+\infty.$\

We recall the following estimate which has been proved in \cite{gp}.

\blem\label{pw}${}$\\
\beqa
\label{pw1}\mrm Pw_\ell&=&w_\ell-\log\l(2\b_\ell^2\d_\ell^{\b_\ell}\r)+4\pi\b_\ell H(\cd,0)+O\l(\d_\ell^{\b_\ell}\r)\\
\nonumber&=&-2\log\l(\d_\ell^{\b_\ell}+|\cd|^{\b_\ell}\r)+4\pi\b_\ell H(\cd,0)+O\l(\d_\ell^{\b_\ell}\r),
\eeqa
and, for any $i,\ell=1,\ds,k$,
\bequ\label{pw2}
\mrm Pw_i(\d_\ell y)=\l\{\bll-2\b_i\log(\d_\ell|y|)+4\pi\b_iH(0,0)+O\l(\fr{1}{|y|^{\b_i}}\l(\fr{\d_i}{\d_\ell}\r)^{\b_i}\r)+O(\d_\ell|y|)+O\l(\d_i^{\b_i}\r)&\tx{if }i<\ell\\-2\b_i\log\d_i-2\log\l(1+|y|^{\b_i}\r)+4\pi\b_iH(0,0)+O(\d_i|y|)+O\l(\d_i^{\b_i}\r)&\tx{if }i=\ell\\-2\b_i\log\d_i+4\pi\b_iH(0,0)+O\l(|y|^{\b_i}\l(\fr{\d_\ell}{\d_i}\r)^{\b_i}\r)+O(\d_\ell|y|)+O\l(\d_i^{\b_i}\r)&\tx{if }i>\ell\earr\r..
\eequ
\elem\

\blem\label{theta}${}$\\
Assume $\b_\ell$ and $\d_\ell$ are defined respectively by \eqref{beta-rec} and \eqref{deltaj}.\\
Then there exists $\g_0>0$ such that, for any $\ell=1,\ds,k$,
\bequ\label{thetao}
|\Th_\ell(y)|=O(\d_\ell|y|+|\la|^{\g_0})\q\q\q\tx{for any }y\in\fr{\mcal A_\ell}{\d_\ell},
\eequ
and in particular
\bequ\label{suptheta}
\sup_{\fr{\mcal A_\ell}{\d_\ell}}|\Th_\ell|=O(1).
\eequ
\elem\

\bpf${}$\\
We will prove the lemma only for odd $\ell,$ i.e. $\ell=2j+1,$ , since the same argument works in the general case. We can also restrict ourselves to consider the case of an odd $k$.\\
We can estimate $\mrm Pw_\ell$ by using Lemma \ref{pw} and then $H$ by the mean value theorem, which gives $H(\d_\ell y,0)=H(0,0)+O(\d_\ell|y|)$:\\
\beqy
\Th_{2j+1}(y)&=&\mrm Pw_{2j+1}(\d_{2j+1}y)-w_{2j+1}(\d_{2j+1}y)+\sum_{i=0}^{j-1}\mrm Pw_{2i+1}(\d_{2j+1}y)+\sum_{i=j+1}^m\mrm Pw_{2i+1}(\d_{2j+1}y)-\fr{a}2\sum_{i=0}^{j-1}\mrm Pw_{2i+2}(\d_{2j+1}y)\\
&-&\fr{a}2\sum_{i=j}^{m-1}\mrm Pw_{2i+1}(\d_{2j+1}y)-(\b_{2j+1}-\a_1)\log|\d_{2j+1}y|-2\pi(\a_1-2)H(\d_{2j+1}y,0)+\log(2\la_1)\\
&=&-\log\l(2\b_{2j+1}^2\r)-\b_{2j+1}\log\d_{2j+1}+4\pi\b_{2j+1}H(0,0)+O(\d_{2j+1}|y|)+O\l(\d_{2j+1}^{\b_{2j+1}}\r)\\
&+&\sum_{i=0}^{j-1}\l(-2\b_{2i+1}\log(\d_{2j+1}|y|)+4\pi\b_{2i+1}H(0,0)+O\l(\fr{1}{|y|^{\b_{2i+1}}}\l(\fr{\d_{2i+1}}{\d_{2j+1}}\r)^{\b_{2i+1}}\r)+O(\d_{2j+1}|y|)+O\l(\d_{2i+1}^{\b_{2i+1}}\r)\r)\\
&+&\sum_{i=j+1}^l\l(-2\b_{2i+1}\log\d_{2i+1}+4\pi\b_{2i+1}H(0,0)+O\l(|y|^{\b_{2i+1}}\l(\fr{\d_{2j+1}}{\d_{2i+1}}\r)^{\b_{2i+1}}\r)+O(\d_{2j+1}|y|)+O\l(\d_{2i+1}^{\b_{2i+1}}\r)\r)\\
&-&\fr{a}2\sum_{i=0}^{j-1}\l(-2\b_{2i+2}\log(\d_{2j+1}|y|)+4\pi\b_{2i+2}H(0,0)+O\l(\fr{1}{|y|^{\b_{2i+2}}}\l(\fr{\d_{2i+2}}{\d_{2j+1}}\r)^{\b_{2i+2}}\r)+O(\d_{2j+1}|y|)+O\l(\d_{2i+2}^{\b_{2i+2}}\r)\r)\\
&-&\fr{a}2\sum_{i=j}^{l-1}\l(-2\b_{2i+2}\log\d_{2i+2}+4\pi\b_{2i+2}H(0,0)+O\l(|y|^{\b_{2i+2}}\l(\fr{\d_{2j+1}}{\d_{2i+2}}\r)^{\b_{2i+2}}\r)+O(\d_{2j+1}|y|)+O\l(\d_{2i+2}^{\b_{2i+2}}\r)\r)\\
&-&2\pi(\a_1-2)H(\d_{2j+1}y,0)-(\b_{2j+1}-\a_1)\log|\d_{2j+1}y|+\log(2\la_1)\\
&=&\ub{-\log\l(2\b_{2j+1}\r)-\b_{2j+1}\log\d_{2j+1}-2\sum_{i=j+1}^m\b_{2i+1}\log\d_{2i+1}+a\sum_{i=j+1}^{m-1}\b_{2i+2}\log\d_{2i+2}}_{=:C_1}\\
&+&\ub{2\pi\l(2\sum_{i=0}^m\b_{2i+1}-a\sum_{i=0}^{m-1}\b_{2i+2}-\a_1+2\r)H(0,0)+\log(2\la_1)}_{=-C_1\tx{ by }\eqref{deltaj}}+\ub{\l(a\sum_{i=0}^{j-1}\b_{2i+1}-2\sum_{i=0}^{j-1}\b_{2i+2}+\a_1-\b_{2j+1}\r)}_{=0\tx{ by }\eqref{beta-rec}}\log(\d_{2j+1}|y|)\\
&+&O(\d_{2j+1}|y|)+\sum_{i=1}^{2m+1}O\l(\d_i^{\b_i}\r)+\sum_{i=1}^{2j}O\l(\fr{1}{|y|^{\b_i}}\l(\fr{\d_i}{\d_{2j+1}}\r)^{\b_i}\r)+\sum_{i=2j+2}^{2m+1}O\l(|y|^{\b_i}\l(\fr{\d_{2j+1}}{\d_i}\r)^{\b_i}\r)\\
&=&O(\d_{2j+1}|y|)+\sum_{i=1}^{2m+1}O\l(\d_i^{\b_i}\r)+\sum_{i=1}^{2j}O\l(\fr{1}{|y|^{\b_i}}\l(\fr{\d_i}{\d_{2j+1}}\r)^{\b_i}\r)+\sum_{i=2j+2}^{2m+1}O\l(|y|^{\b_i}\l(\fr{\d_{2j+1}}{\d_i}\r)^{\b_i}\r)\\
&=&O(\d_{2j+1}|y|)+\sum_{i=1}^{2m+1}O\l(\d_i^{\b_i}\r)+\sum_{i=1}^{2j}O\l(\l(\fr{\d_i^2}{\d_{2j}\d_{2j+1}}\r)^\fr{\b_i}2\r)+\sum_{i=2j+2}^{2m+1}O\l(\l(\fr{\d_{2j+1}\d_{2j+2}}{\d_i^2}\r)^\fr{\b_i}2\r)\\
&=&O(\d_{2j+1}|y|)+\sum_{i=1}^{2m+1}O\l(\d_i^{\b_i}\r)+\sum_{i=1}^{2j}O\l(\l(\fr{\d_{2j}}{\d_{2j+1}}\r)^\fr{\b_i}2\r)+\sum_{i=2j+2}^{2m+1}O\l(\l(\fr{\d_{2j+1}}{\d_{2j+2}}\r)^\fr{\b_i}2\r)\\
&=&O(\d_{2j+1}|y|)+O\l(\min_i\d_i^{\b_i}\r)+O\l(\min_{i,\ell}\l(\fr{\d_\ell}{\d_{\ell+1}}\r)^\fr{\b_i}2\r)\\
&=&O(\d_{2j+1}|y|)+O(|\la|^{\g_0})\\
&=&O(\d_{2j+1}|y|+|\la|^{\g_0}),
\eeqy
where we used that $\sqrt\fr{\d_{2j}}{\d_{2j+1}}\le|y|\le\sqrt\fr{\d_{2j+1}}{\d_{2j+2}}$ and that $\d_i\le\d_{2j}$ and $\d_{2j+2}\le\d_{i'}$ for any $i<2j+1<i'$.\\
\eqref{suptheta} follows straightforwardly from \eqref{thetao}, since $\d_\ell|y|=O(1)$ for any $y\in\fr{\mcal A_\ell}{\d_\ell}$.
\epf\

\subsection{Chebyshev polynomials and the $\beta_\ell$'s}
\label{polinomi}

In this sub-section we shall prove that the $\beta_\ell$'s defined in \eqref{betaj} and in \eqref{beta-rec} coincide.

Let us introduce the polynomials
\bequ\label{pij}
\left\{\begin{aligned}
 P_1(t)&=1\\
 P_2(t)&=1\\
&\ \vdots \\
 P_{2j+1}(t)&= \sum_{i=0}^j(-1)^{j+i}{j+i\choose 2i}t^i\\
P_{2j+2}(t)&= \sum_{i=0}^j(-1)^{j+i}{j+i+1\choose 2i+1}t^i\\
&\ \vdots\end{aligned}\right.
\eequ\

By induction, is not difficult to check that the real numbers defined in \eqref{beta-rec} satisfy \eqref{betaj}, since

\bequ\label{bj}
\left\{\begin{aligned}
\b_1&=\a_1\\
\b_2&=b\a_1+\a_2\\
& \ \vdots\\
\b_{2j+1} &= \a_1\sum_{i=0}^j(-1)^{j+i}{j+i\choose 2i}a^ib^i+\a_2\sum_{i=0}^{j-1}(-1)^{j+i-1}{j+i\choose 2i+1}a^{i+1}b^i=\a_1P_{2j+1}(ab)+a\a_2P_{2j}(ab)\\
\b_{2j+2}&= \a_1\sum_{i=0}^j(-1)^{j+i}{j+i+1\choose 2i+1}a^ib^{i+1}+\a_2\sum_{i=0}^j(-1)^{j+i}{j+i\choose 2i}a^ib^i=b\a_1P_{2j+2}(ab)+\a_2P_{2j+1}(ab).\\
&\ \vdots
\end{aligned}\right.
\eequ\

Therefore, the problem reduces to prove that the polynomials defined in \eqref{pij} coincide with the polynomial defined in \eqref{pj}.\

Now, the polynomials defined in \eqref{pij} can be expressed in terms of Chebyshev's polynomials 
\bequ\label{tj}
\l\{\bll
T_0(x)=1,\\
T_1(x)=x,\\
T_{\ell+1}(x)=2xT_\ell(x)-T_{\ell-1}(x)& \hbox{if $\ell\ge2$}
\earr\r..
\eequ\

\begin{lemma}${}$\\
Let $P_\ell$ be defined by \eqref{pj} and $T_\ell$ be defined by \eqref{tj}.\\
Then, for any $j\in\N$ and $x\in\R$ it holds:
\bequ\label{tjpj}T_{2j+1}(x)=1+(x-1)\l(P_{2j+1}(2x+2)\r)^2\q\q \hbox{and}\q\q T_{2j+2}(x)=1+\l(2x^2-2\r)\l(P_{2j+2}(2x+2)\r)^2.\eequ
\end{lemma}\

\bpf${}$\\
We proceed by induction. We can easily see that the Proposition is true for $\ell=1,2$.\\
Let us now assume the Proposition to hold for any positive integer up to $2j$ and let us show it still holds true for $2j+1$ and $2j+2$.\\
First of all, by induction we can easily show that $P_\ell$ verifies the following properties:
\bequ
\label{pj1}P_{2j+1}(t)=tP_{2j}(t)-P_{2j-1}(t)\q\q\q P_{2j+2}(t)=P_{2j+1}(t)-P_{2j}(t),\eequ
and also
\bequ
\label{pj2}\l(P_{2j+1}(t)\r)^2+t\l(P_{2j}(t)\r)^2-tP_{2j+1}(t)P_{2j}(t)-1=0 \q\q\q t\l(P_{2j+2}(t)\r)^2+\l(P_{2j+1}(t)\r)^2-tP_{2j+2}(t)P_{2j+1}(t)-1=0.
\eequ

Using \eqref{pj1} and \eqref{pj2} we get, for odd indexes:
\beqy
&&T_{2j+1}(x)\\
&=&2xT_{2j}(x)-T_{2j-1}(x)\\
&=&2x+\l(4x^3-4x\r)P_{2j}(2x+2)^2-1-(x-1)P_{2j-1}(2x+2)^2\\
&=&1+(x-1)\l(\l(4x^2+4x\r)P_{2j}(2x+2)^2-P_{2j-1}(2x+2)^2+2\r)\\
&=&1+(x-1)\l(P_{2j+1}(2x+2)^2-2\l((2x+2)P_{2j}(2x+2)^2-P_{2j-1}(2x+2)^2-(2x+2)P_{2j}(2x+2)P_{2j-1}(2x+2)-1\r)\r)\\
&=&1+(x-1)P_{2j+1}(2x+2)^2;
\eeqy
similarly, for even indexes:
\beqy
&&T_{2j+2}(x)\\
&=&2xT_{2j+1}(x)-T_{2j}(x)\\
&=&2x+\l(2x^2-2x\r)P_{2j+1}(2x+2)^2-1-\l(2x^2-2\r)P_{2j}(2x+2)^2\\
&=&1+(2x-2)\l(xP_{2j+1}(2x+2)^2-(x+1)P_{2j}(2x+2)^2+1\r)\\
&=&1+(2x-2)\l((x+1)P_{2j+2}(2x+2)^2-\l(P_{2j+1}(2x+2)^2+(2x+2)P_{2j}(2x+2)^2-(2x+2)P_{2j+2}(2x+2)P_{2j}(2x+2)-1\r)\r)\\
&=&1+\l(2x^2-2\r)P_{2j+2}(2x+2)^2.
\eeqy
\epf\

\begin{remark}\label{tmax2}${}$\\
By using the properties of Chebyshev's polynomials (see for instance \cite{riv}), we easily find 
the explicit expression of $T_\ell$ as
$$T_\ell(x)=1+2^{\ell-1}\prod_{i=1}^\ell\l(x-\cos\fr{2\pi i}{\ell}\r).$$
which can be rewritten, if $\ell=2j+1$ is odd or $\ell=2j+2$ is even, as
\bequ\label{tjpj2}T_{2j+1}(x)=1+2^{2j}(x-1)\prod_{i=1}^j\l(x-\cos\fr{2\pi i}{2j+1}\r)^2\q \hbox{or}\q T_{2j+2}(x)=1+2^{2j+1}\l(x^2-1\r)\prod_{i=1}^j\l(x-\cos\fr{\pi i}{j+1}\r)^2,\eequ
respectively.\\
Now, if we compare \eqref{tjpj} and \eqref{tjpj2} we get the explicit expression for $P_\ell$ given in \eqref{pj}.
\end{remark}\

\begin{remark}\label{tmax3}${}$\\
If $ab=2\cos\l(\fr{2\pi}k\r)+2$ for some $k\in\N$, then we have $P_\ell(ab)>0$ for any $\ell=1,\ds,k-1$ and $P_k(ab)=0>P_{k+1}(ab)$; hence, by the definitions \eqref{betaj} of $\b_\ell$ and \eqref{kmax} of $k_{\max}$ we get $k_{\max}=k=\fr{2\pi}{\arccos\l(\fr{ab}2-1\r)}$.\\
On the other hand, if $2\cos\l(\fr{2\pi}k\r)+2<ab<2\cos\l(\fr{2\pi}{k+1}\r)+2$, then $P_\ell(ab)>0$ for $\ell=1,\ds,k-1$ and $P_k(ab),P_{k+1}(ab)<0$; hence, $\b_\ell>0$ for $\ell\le k-1$ and $\b_{k+1}<0$, so $k_{\max}$ could be either $k-1$ or $k$.\\
Finally, if $ab\ge4$, then clearly $P_\ell(ab)>0$ for all $\ell$, hence $b_\ell>0$ and $k_{\max}=+\infty$.\\
\end{remark}\

\begin{remark}\label{t=4}${}$\\
By \eqref{tjpj} we immediately deduce that
$$\l(P_{2j+1}(4)\r)^2=T'_{2j+1}(1)=(2j+1)^2\q \hbox{and}\q \l(P_{2j+2}(4)\r)^2=\frac 1 4T'_{2j+2}(1)=\frac{(2j+2)^2}4=(j+1)^2,$$
because
 the
 Chebyshev's polynomials satisfy (by induction, for instance)
 $T'_\ell(1)=\ell ^2$ for any $\ell\ge1.$
\end{remark}\

\begin{remark}\label{masse}
${}$\\
The validity of \eqref{mass1} and \eqref{mass2} follows by the fact 
that the coefficients $\b_\ell$ verify the following properties (by induction, for instance):
\begin{eqnarray*}
 \sum_{i=0}^{2j+1}\b_{2i+1}=P_{2j+1}(ab)(\a_1P_{2j+1}(ab)+a\a_2P_{2j}(ab)),&\q&\sum_{i=0}^{2j+2}\b_{2i+1}=aP_{2j+2}(ab)(b\a_1P_{2j+2}(ab)+\a_2P_{2j+1}(ab)),\\
\nonumber\sum_{i=0}^{2j+1}\b_{2i+2}=P_{2j+1}(ab)(b\a_1P_{2j+2}(ab)+\a_2P_{2j+1}(ab)),&\q&\sum_{i=0}^{2j+2}\b_{2i+2}=bP_{2j+2}(ab)(\a_1P_{2j+3}(ab)+a\a_2P_{2j+2}(ab)). 
\end{eqnarray*}
\end{remark}\

\section{Proof of the main theorem}\label{proof}\

In this section, we prove the existence of a solution to system \eqref{plambda} using a contraction mapping argument and we study its properties.

\bprop\label{contr}${}$\\
There exist $\g,R,\ol\la>0$ such that for any $\la\in\l(0,\ol\la\r)\x\l(0,\ol\la\r)$ there exists a unique $\phi_\la=(\phi_{1,\la},\phi_{2,\la})\in\mbf H$ such that:
\bite
\item $W_\la+\phi_\la$ solves \eqref{plambda}, namely
$$\l\{\bll-\D(W_{1,\la}+\phi_{1,\la})=2\la_1h_1e^{W_{1,\la}+\phi_{1,\la}}-a\la_2h_2e^{W_{2,\la}+\phi_{2,\la}}&\tx{in }\O\\
-\D(W_{2,\la}+\phi_{2,\la})=2\la_2h_2e^{W_{2,\la}+\phi_{2,\la}}-b\la_1h_1e^{W_{1,\la}+\phi_{1,\la}}&\tx{in }\O\earr\r.;$$
\item $\|\phi_\la\|\le R|\la|^\g\log\fr{1}{|\la|}$.
\eite
\eprop\

\bpf${}$\\
We point out that $W_\la+\phi_\la$ solves \eqref{plambda} if and only if
$$\mcal L_\la\phi=\mcal N_\la(\phi)-\mcal S_\la\phi-\mcal R_\la.$$

where the linear operator $\mcal L_\la: H^1_0(\O)\x H^1_0(\O)\to L^p(\O)\x L^p(\O)$ is defined by
\bequ\label{llambda}
\mcal L_\la(\phi):=\l(\barr{c}\mcal L_{1,\la}\phi\\\mcal L_{2,\la}\phi\earr\r)=\l(\barr{c}-\D\phi_1-\l(\sum_{j=0}^{\l[\fr{k-1}2\r]}2\b_{2j+1}^2\fr{\d_{2j+1}^{\b_{2j+1}}|\cd|^{\b_{2j+1}-2}}{\l(\d_{2j+1}^{\b_{2j+1}}+|\cd|^{\b_{2j+1}}\r)^2}\phi_1-\fr{a}2\sum_{j=0}^{\l[\fr{k-2}2\r]}2\b_{2j+2}^2\fr{\d_{2j+2}^{\b_{2j+2}}|\cd|^{\b_{2j+2}-2}}{\l(\d_{2j+2}^{\b_{2j+2}}+|\cd|^{\b_{2j+2}}\r)^2}\phi_2\r)\\-\D\phi_2-\l(\sum_{j=0}^{\l[\fr{k-2}2\r]}2\b_{2j+1}^2\fr{\d_{2j+2}^{\b_{2j+2}}|\cd|^{\b_{2j+2}-2}}{\l(\d_{2j+2}^{\b_{2j+2}}+|\cd|^{\b_{2j+2}}\r)^2}\phi_2-\fr{b}2\sum_{j=0}^{\l[\fr{k-1}2\r]}2\b_{2j+1}^2\fr{\d_{2j+1}^{\b_{2j+1}}|\cd|^{\b_{2j+1}-2}}{\l(\d_{2j+1}^{\b_{2j+1}}+|\cd|^{\b_{2j+1}}\r)^2}\phi_1\r)\earr\r),
\eequ
the error function $\mcal R_\la\in L^p(\O)\x L^p(\O)$ is defined by
\bequ\label{rlambda}
\mcal R_\la:=\l(\barr{c}\mcal R_{1,\la}\\\mcal R_{2,\la}\earr\r)=\l(\barr{c}-\D W_{1,\la}-2\la_1h_1e^{W_{1,\la}}+a\la_2h_2e^{W_{2,\la}}\\-\D W_{2,\la}-2\la_2h_2e^{W_{2,\la}}+b\la_1h_1e^{W_{1,\la}}\earr\r),
\eequ
the error linear operator $\mcal S_\la:H^1_0(\O)\x H^1_0(\O)\to L^p(\O)\x L^p (\O)$ is defined by
\bequ\label{slambda}
\mcal S_\la\l(\barr{c}\phi_1\\\phi_2\earr\r):=\l(\barr{c}\l(\sum_{j=0}^{\l[\fr{k-1}2\r]}|\cd|^{\b_{2j+1}}e^{w_{2j+1}}-2\la_1h_1e^{W_{1,\la}}\r)\phi_1-\fr{a}2\l(\sum_{j=0}^{\l[\fr{k-2}2\r]}|\cd|^{\b_{2j+2}}e^{w_{2j+2}}-2\la_2h_2e^{W_{2,\la}}\r)\phi_2\\\l(\sum_{j=0}^{\l[\fr{k-2}2\r]}|\cd|^{\b_{2j+2}}e^{w_{2j+2}}-2\la_2h_2e^{W_{2,\la}}\r)\phi_2-\fr{b}2\l(\sum_{j=0}^{\l[\fr{k-1}2\r]}|\cd|^{\b_{2j+1}}e^{w_{2j+1}}-2\la_1h_1e^{W_{1,\la}}\r)\phi_1\earr\r) 
\eequ
and the quadratic term
 $\mcal N_\la :H_0^1(\O)\x H^1_0(\O)\to L^p(\O)\x L^p(\O)$ is defined by
\bequ\label{nlambda}
\mcal N_\la(\phi):=\l(\barr{c}\mcal N_{1,\la}(\phi)\\\mcal N_{2,\la}(\phi)\earr\r)=\l(\barr{c}2\la_1h_1e^{W_{1,\la}}\l(e^{\phi_1}-1-\phi_1\r)-a\la_2h_2e^{W_{2,\la}}\l(e^{\phi_2}-1-\phi_2\r)\\2\la_2h_2e^{W_{2,\la}}\l(e^{\phi_2}-1-\phi_2\r)-b\la_1h_1e^{W_{1,\la}}\l(e^{\phi_1}-1-\phi_1\r)\earr\r).
\eequ

Since Proposition \ref{lla} ensures that $\mcal L_\la:\mbf H\to\mbf H$ is invertible, this is equivalent to requiring $\phi_\la$ to be a fixed point of the map
$$\mcal T_\la:\phi\mapsto(\mcal L_\la)^{-1}(\mcal N_\la(\phi)-\mcal S_\la\phi-\mcal R_\la);$$
therefore, the existence of such a $\phi_\la$ will follow by showing that $\mcal T_\la$ is a contraction on the ball
$$B_{\g,\la,R}:=\l\{\phi\in\mbf H:\,\|\phi\|\le R|\la|^\g\log\fr{1}{|\la|}\r\},$$
for $\g,\la$ small enough and $R$ large enough.\\

We first show that $\mcal T_\la$ maps $B_{\g,\la,R}$ into itself.\\
We will use Proposition \ref{lla} to get estimates on $\mcal L_\la$ and Lemmas \ref{rla}, \ref{sla}, \ref{nla} to estimate $\mcal N_\la(\phi),\,\mcal S_\la\phi,\,\mcal R_\la$, respectively.\\
With the notation of these Lemmas, we take $\g\le\max\{\g_1,\g_2\}$ and $p$ so close to $1$ that all such Lemmas apply and $\g_3(1-p)+\g>0$. We then take $C$ as in Lemma \ref{nla} and $\ol\la$ so small that $e^CR^2\ol\la^{\g_3(1-p)+\g}\l(\log\fr{1}{\ol\la}\r)^2\le1$. Finally, we take $R>0$ greater than the three constants which define the $O$ in Lemmas \ref{rla}, \ref{sla}, \ref{nla}, times the $C$ appearing in Proposition \ref{lla}.\\
Notice that these choice imply that $R|\la|^\g\log\fr{1}{|\la|}\le1$; therefore
\beqy
\|\mcal T_\la(\phi)\|&\le&C\log\fr{1}{|\la|}\l(\l\|\mcal N_\la(\phi)\|_p+\|\mcal S_\la\phi\|_p+\|\mcal R_\la\r\|_p\r)\\
&\le&C\log\fr{1}{|\la|}\l(|\la|^{\g_3(1-p)}\|\phi\|^2e^{C\|\phi\|^2}+|\la|^{\g_2}\|\phi\|+|\la|^{\g_1}\r)\\
&\le&C\log\fr{1}{|\la|}\l(R^2e^C\l(\log\fr{1}{|\la|}\r)^2|\la|^{\g_3(1-p)+2\g}+|\la|^\g\r)\\
&\le&C|\la|^\g\log\fr{1}{|\la|}\\
&\le&R|\la|^\g\log\fr{1}{|\la|}.
\eeqy\

Moreover, we also get
\beqy
\|\mcal T_\la(\phi)-\mcal T_\la(\phi')\|&\le&C\log\fr{1}{|\la|}\l(\l\|\mcal N_\la(\phi)-\mcal N_\la(\phi')\|_p+\|\mcal S_\la(\phi-\phi')\r\|_p\r)\\
&\le&C\log\fr{1}{|\la|}\l(|\la|^{\g_3(1-p)}\|\phi-\phi'\|(\|\phi\|+\|\phi'\|)e^{C\l(\|\phi\|^2+\l\|\phi\r\|^2\r)}+|\la|^{\g_2}\|\phi-\phi'\|\r)\\
&\le&C\l(2Re^{2C}|\la|^{\g_3(1-p)+\g}\l(\log\fr{1}{|\la|}\r)^2+|\la|^{\g_2}\log\fr{1}{|\la|}\r)\|\phi-\phi'\|\\
&\le&C\l(2\fr{e^C}R+\ol\la^{\g_2}\log\fr{1}{\ol\la}\r)\|\phi-\phi'\|;
\eeqy
with the constant multiplying $\|\phi-\phi'\|$ being smaller than $1$, after taking larger $R$ and/or smaller $\ol\la$, if needed. This concludes the proof.
\epf\

\bpf[Proof of Theorem \ref{main}]${}$\\
By Proposition \ref{contr} we get $u_\la=W_\la+\phi_\la$ which solves \eqref{plambda}.\\

Let us prove \eqref{mass}.\\
We basically show that $\phi_\la$ is negligible in this computations, thanks to the estimates from Proposition \ref{nla} and in particular \eqref{lp}. Then, we compare $W_{i,\la}$ with $w_{2j+i}$ using the estimate \eqref{w} from Lemma \ref{sla}
\beqy
&&\l|\la_i\int_{B_r(0)}h_ie^{u_{i,\la}}-2\pi\sum_{j=0}^{\l[\fr{k-i}2\r]}\b_{2j+i}\r|\\
&=&\l|\la_i\int_{B_r(0)}h_ie^{W_{i,\la}+\phi_{i,\la}}-\fr{1}2\sum_{j=0}^{\l[\fr{k-i}2\r]}\int_{\R^2}2\b_{2j+i}^2\fr{|\cd|^{\b_{2j+i}-2}}{\l(1+|\cd|^{\b_{2j+i}}\r)^2}\r|\\
&\le&\int_\O\la_ih_ie^{W_{i,\la}}\l|e^{\phi_{i,\la}}-1\r|+\l|\la_i\int_{B_r(0)}h_ie^{W_{i,\la}}-\fr{1}2\sum_{j=0}^{\l[\fr{k-i}2\r]}\int_{\R^2}2\b_{2j+i}^2\fr{|\cd|^{\b_{2j+i}-2}}{\l(1+|\cd|^{\b_{2j+i}}\r)^2}\r|\\
&\le&\int_\O\la_ih_ie^{W_{i,\la}}|\phi_{i,\la}|e^{\phi_{i,\la}}+\l|\la_i\int_{B_r(0)}h_ie^{W_{i,\la}}-\fr{1}2\sum_{j=0}^{\l[\fr{k-i}2\r]}\int_{B_\fr{r}{\d_{2j+i}}(0)}2\b_{2j+i}^2\fr{|\cd|^{\b_{2j+i}-2}}{\l(1+|\cd|^{\b_{2j+i}}\r)^2}\r|+o(1)\\
&\le&\l\|\la_ih_ie^{W_{i,\la}}\r\|_p\|\phi_{i,\la}\|_q\l\|e^{\phi_{i,\la}}\r\|_\fr{pq}{pq-p-q}+\int_{B_r(0)}\l|\la_ih_ie^{W_{i,\la}}-\fr{1}2\sum_{j=0}^{\l[\fr{k-i}2\r]}|\cd|^{\b_{2j+i}-2}e^{w_{2j+i}}\r|+o(1)\\
&\le&CR|\la|^{\g_3(1-p)+\g}\log\fr{1}{|\la|}+C|\la|^{\g_1}+o(1)\\
&\us{\la\to0}\to&0.
\eeqy
Since this holds true for any $r$, then letting $r$ tend to $0$ we find the value of $m_i(0)$.\\

Finally, we prove \eqref{green}.\\
First of all, $u_\la$ is bounded in $W^{1,q}(\O)\x W^{1,q}(\O)$ for any $q<2$, because $W^{1,\fr{q}{q-1}}(\O)\inc C\l(\ol\O\r)$, hence for any $\ph\in W^{1,\fr{q}{q-1}}(\O)$ with $\|\ph\|_{W^{1,\fr{q}{q-1}}(\O)}\le1$ we have
$$\l|\int_\O\n u_{i,\la}\cd\n\ph\r|=\l|\int_\O(-\D u_{i,\la})\ph\r|\le C\l(\la_1\int_\O h_1e^{u_{1,\la}}+\la_2\int_\O h_2e^{u_{2,\la}}\r)\|\ph\|_\infty\le C.$$
From \eqref{pw1} we get $\mrm Pw_\ell(x)\us{\la\to0}\to4\pi\b_\ell G(\cd,0)$ pointwise in $\O\sm\{0\}$. Since $\l\|\phi_\la\r\|\us{\la\to0}\to0$, from the definition of $u_\la$ and $m_1,m_2$ we deduce that the weak limit of $u_\la$ in $W^{1,q}(\O)$ must be the one in \eqref{green}.\\
Moreover, from \eqref{pw2} and the definition of $W_{i,\la}$ we deduce that the latter are both bounded in $L^\infty_{\mrm{loc}}(\O\sm\{0\})$. Therefore, for any $\mcal K\Sub\O\sm\{0\}$,
$$\int_{\mcal K}|-\D u_{i,\la}|^q\le C\l(\int_{\mcal K}|\cd|^{q(\a_1-2)}e^{q\l(W_{1,\la}+\phi_{1,\la}\r)}+\int_{\mcal K}|\cd|^{q(\a_2-2)}e^{q\l(W_{2,\la}+\phi_{2,\la}\r)}\r)\le Ce^{q\|W_\la\|_\infty}\l(\int_{\mcal K}e^{q\phi_{1,\la}}+\int_{\mcal K}e^{q\phi_{2,\la}}\r)\le C.$$
Therefore, a standard bootstrap method will imply convergence in $C^\infty(\mcal K)$ hence, being $\mcal K$ arbitrary, in $C^\infty_{\mrm{loc}}(\O\sm\{0\})$.
\epf\

\section{The error terms}\label{errors}\

In this section we estimate in the $L^p$ norm the function $\mcal R_\la$ defined in \eqref{rlambda}, the linear operator $\mcal S_\la$ defined in \eqref{slambda} and the quadratic term $\mcal N_\la$ defined in \eqref{nlambda}.\\
Roughly speaking, both $\mcal R_\la$ and $\mcal S_\la$ will decay as a power of $\la$ if $p$ is close enough to $1$. On the other hand, the norm of $\mcal N_\la$ will diverge as $\la$ goes to $0$, but its growth will be slow for small $p$.\\
The estimates for $\mcal S_\la$ and $\mcal N_\la$ will require mostly the same calculations as the ones needed for $\mcal R_\la$.

\subsection{The function $\mcal R_\la$}

\blem\label{rla}${}$\\
Let $\mcal R_\la$ be defined by \eqref{rlambda}.\\
There exists $p_0>1$ and $\g_1>0$ such that for any $p\in[1,p_0)$
$$\|\mcal R_\la\|_p=O(|\la|^{\g_1}).$$
\elem\

\bpf${}$\\
We will only provide estimates for $\mcal R_{1,\la}$, since the estimates for $\mcal R_{2,\la}$ are similar.\\
First of all, by the very definition of $W_{i,\la}$ and triangular inequalities, we can split the $L^p$ norm of $\mcal R_{1,\la}$ in the following way:
\beqy
\int_\O\l|\mcal R_{1,\la}\r|^p&=&\int_\O\l|\sum_{j=0}^{\l[\fr{k-1}2\r]}|\cd|^{\b_{2j+1}-2}e^{w_{2j+1}}-\fr{a}2\sum_{j=0}^{\l[\fr{k-2}2\r]}|\cd|^{\b_{2j+2}-2}e^{w_{2j+2}}-2\la_1h_1e^{\sum_{m=0}^{\l[\fr{k-1}2\r]}\mrm Pw_{2m+1}-\fr{a}2\sum_{m=0}^{\l[\fr{k-2}2\r]}\mrm Pw_{2m+2}}\r.\\
&+&\l.a\la_2h_2e^{\sum_{m=0}^{\l[\fr{k-2}2\r]}\mrm Pw_{2m+2}-\fr{b}2\sum_{m=0}^{\l[\fr{k-1}2\r]}\mrm Pw_{2m+1}}\r|^p\\
&\le&C\int_\O\l|\sum_{j=0}^{\l[\fr{k-1}2\r]}|\cd|^{\b_{2j+1}-2}e^{w_{2j+1}}-2\la_1h_1e^{\sum_{m=0}^{\l[\fr{k-1}2\r]}\mrm Pw_{2m+1}-\fr{a}2\sum_{m=0}^{\l[\fr{k-2}2\r]}\mrm Pw_{2j+2}}\r|^p\\
&+&C\int_\O\l|\sum_{j=0}^{\l[\fr{k-2}2\r]}|\cd|^{\b_{2j+2}-2}e^{w_{2j+2}}-2\la_2h_2e^{\sum_{m=0}^{\l[\fr{k-2}2\r]}\mrm Pw_{2m+2}-\fr{b}2\sum_{m=0}^{\l[\fr{k-1}2\r]}\mrm Pw_{2m+1}}\r|^p\\
&\le&C\sum_{j=0}^{\l[\fr{k-1}2\r]}\ub{\int_{A_{2j+1}}\l||\cd|^{\b_{2j+1}-2}e^{w_{2j+1}}-2\la_1h_1e^{\sum_{m=0}^{\l[\fr{k-1}2\r]}\mrm Pw_{2m+1}-\fr{a}2\sum_{m=0}^{\l[\fr{k-2}2\r]}\mrm Pw_{2m+2}}\r|^p}_{=:I'_{2j+1}}\\
&+&C\sum_{j=0}^{\l[\fr{k-1}2\r]}\sum_{i=1,i\ne 2j+1}^k\ub{\int_{\mcal A_i}\l||\cd|^{\b_{2j+1}-2}e^{w_{2j+1}}\r|^p}_{=:I''_{i,2j+1}}+C\sum_{i=0}^{\l[\fr{k-2}2\r]}\ub{\int_{\mcal A_{2i+2}}\l|\la_1h_1e^{\sum_{m=0}^{\l[\fr{k-1}2\r]}\mrm Pw_{2m+1}-\fr{a}2\sum_{m=0}^{\l[\fr{k-2}2\r]}\mrm Pw_{2m+2}}\r|^p}_{=:I'''_{2i+2}}\\
&+&C\sum_{j=0}^{\l[\fr{k-2}2\r]}\ub{\int_{\mcal A_{2j+2}}\l||\cd|^{\b_{2j+2}-2}e^{w_{2j+2}}-2\la_2h_2e^{\sum_{m=0}^{\l[\fr{k-2}2\r]}\mrm Pw_{2m+2}-\fr{b}2\sum_{m=0}^{\l[\fr{k-1}2\r]}\mrm Pw_{2m+1}}\r|^p}_{=:I'_{2j+2}}\\
&+&C\sum_{j=0}^{\l[\fr{k-2}2\r]}\sum_{i=1,i\ne 2j+2}^k\ub{\int_{\mcal A_i}\l||\cd|^{\b_{2j+2}-2}e^{w_{2j+2}}\r|^p}_{=:I''_{i,2j+2}}+C\sum_{i=0}^{\l[\fr{k-1}2\r]}\ub{\int_{\mcal A_{2i+1}}\l|\la_2h_2e^{\sum_{m=0}^{\l[\fr{k-2}2\r]}\mrm Pw_{2m+2}-\fr{b}2\sum_{m=0}^{\l[\fr{k-1}2\r]}\mrm Pw_{2m+1}}\r|^p}_{=:I'''_{2i+1}}\\.
\eeqy
Now we suffice to estimate separately each of $I'_\ell,I''_{i,\ell},I'''_{i}$.\\
To handle with $I'_\ell$ we use the definition \eqref{thetaj} of $\Th_\ell$ and Lemma \ref{theta}:
\beqy
I'_\ell&=&\int_{\mcal A_\ell}\l||x|^{\b_\ell-2}e^{w_\ell(x)}\l(1-e^{\Th_\ell\l(\fr{x}{\d_\ell}\r)}\r)\r|^p\mrm dx\\
&=&\l(2\b_\ell^2\r)^p\d_\ell^{2-2p}\int_{\fr{\mcal A_\ell}{\d_\ell}}\fr{|y|^{(\b_\ell-2)p}}{(1+|y|^{\b_\ell})^{2p}}\l|1-e^{\Th_\ell(y)}\r|^p\mrm dy\\
&\le&\l(2\b_\ell^2\r)^p\d_\ell^{2-2p}\int_{\fr{\mcal A_\ell}{\d_\ell}}\fr{|y|^{(\b_\ell-2)p}}{(1+|y|^{\b_\ell})^{2p}}|\Th_\ell(y)|^pe^{p|\Th_\ell(y)|}\mrm dy\\
&\le&C\d_\ell^{2-2p}\int_{\fr{\mcal A_\ell}{\d_\ell}}\fr{|y|^{(\b_\ell-2)p}}{(1+|y|^{\b_\ell})^{2p}}|\Th_\ell(y)|^p\mrm dy\\
&\le&C\d_\ell^{2-2p}\int_{\fr{\mcal A_\ell}{\d_\ell}}\fr{|y|^{(\b_\ell-2)p}}{(1+|y|^{\b_\ell})^{2p}}\l|\d_\ell y+|\la|^{\g_0}\r|^p\mrm dy\\
&\le&C\l(\d_\ell^{2-\min\l\{1,2-\fr{\b_\ell}2\r\}p}\int_{\fr{\mcal A_\ell}{\d_\ell}}\fr{|y|^{\max\l\{\fr{3}2\b_\ell-2,\b_\ell-1\r\}p}}{(1+|y|^{\b_\ell})^{2p}}\mrm dy+\d_\ell^{2-2p}|\la|^{p\g_0}\int_{\fr{\mcal A_\ell}{\d_\ell}}\fr{|y|^{(\b_\ell-2)p}}{(1+|y|^{\b_\ell})^{2p}}\mrm dy\r)\\
&\le&C\l(\d_\ell^{2-\min\l\{1,2-\fr{\b_\ell}2\r\}p}+\d_\ell^{2-2p}|\la|^{p\g_0}\r),
\eeqy
which can be estimated by a power of $|\la|$ if $p$ is close enough to $1$.\\
Concerning $I''_{i,\ell}$, we have:
\beqa
\nonumber I''_{i,\ell}&=&\int_{\mcal A_i}\l|\fr{2\b_\ell^2\d_\ell^{\b_\ell}|x|^{\b_\ell-2}}{\l(\d_\ell^{\b_\ell}+|x|^{\b_\ell}\r)^2}\r|^p\mrm dx\\
\nonumber&=&\l(2\b_\ell^2\r)^p\d_\ell^{2-2p}\int_{B_{\sqrt\fr{\d_i\d_{i+1}}{\d_\ell}}\sm B_{\sqrt\fr{\d_{i-1}\d_i}{\d_\ell}}}\fr{|y|^{(\b_\ell-2)p}}{(1+|y|^{\b_\ell})^{2p}}\mrm dy\\
\nonumber&\le&C\d_\ell^{2-2p}\l\{\bll\l(\fr{\sqrt{\d_i\d_{i+1}}}{\d_\ell}\r)^{(\b_\ell-2)p+2}&\tx{if }i<\ell\\\l(\fr{\d_\ell}{\sqrt{\d_{i-1}\d_i}}\r)^{-(\b_\ell+2)p+2}&\tx{if }\ell>i\earr\r.\\
\nonumber&\le&C\d_\ell^{2-2p}\l\{\bll\l(\fr{\d_{\ell-1}}{\d_\ell}\r)^\fr{(\b_\ell-2)p+2}2&\tx{if }i<\ell\\\l(\fr{\d_\ell}{\d_{\ell+1}}\r)^\fr{-(\b_\ell+2)p+2}2&\tx{if }\ell>i\earr\r.\\
\label{i2}&\le&C\d_\ell^{2-2p}|\la|^{\g'},
\eeqa
which is still bounded by a power of $\la$ for small $p$.\\
Finally, for $I'''_{2i+2}$, we use \eqref{pw1}, the fact that $\d_\ell\le|y|\le\d_{\ell'}$ for any $y\in\mcal A_{2i+2}$ and $\ell<2i+2<\ell'$, and then the properties \eqref{deltaj} of $\d_\ell$'s and \eqref{bj} of $\b_\ell$:
\beqy
I'''_{2i+2}&\le&C\int_{\mcal A_{2i+2}}\l|\la_1|x|^{\a_1-2}\prod_{m=0}^{\l[\fr{k-1}2\r]}\fr{1}{\l(\d_{2m+1}^{\b_{2m+1}}+|x|^{\b_{2m+1}}\r)^2}\prod_{m=0}^{\l[\fr{k-2}2\r]}\l(\d_{2m+2}^{\b_{2m+2}}+|x|^{\b_{2m+2}}\r)^a\r|^p\mrm dx\\
&\le&C\la_1^p\prod_{m=i+1}^{\l[\fr{k-1}2\r]}\d_{2m+1}^{-2\b_{2m+1}p}\prod_{m=i+1}^{\l[\fr{k-2}2\r]}\d_{2m+2}^{a\b_{2m+2}p}\int_{\mcal A_{2i+2}}|x|^{(\a_1-2)p-2\sum_{m=0}^i\b_{2m+1}p+a\sum_{m=0}^{i-1}\b_{2m+2}p}\l(\d_{2i+2}^{\b_{2i+2}}+|x|^{\b_{2i+2}}\r)^{ap}\mrm dx\\
&=&C\d_{2i+3}^{-\b_{2i+3}p}\int_{\mcal A_{2i+2}}|x|^{-(\b_{2i+1}+2)p}\l(\d_{2i+2}^{\b_{2i+2}}+|x|^{\b_{2i+2}}\r)^{ap}\mrm dx\\
&=&C\d_{2i+3}^{-\b_{2i+3}p}\d_{2i+2}^{2+(-2-\b_{2i+1}+a\b_{2i+2})p}\int_{\fr{\mcal A_{2i+2}}{\d_{2i+2}}}|y|^{-(\b_{2i+1}+2)p}\l(1+|y|^{\b_{2i+2}}\r)^{ap}\mrm dy\\
&=&C\l(\fr{\d_{2i+2}}{\d_{2i+3}}\r)^{\b_{2i+3}p}\d_{2i+2}^{2-2p}\int_{B_{\sqrt\fr{\d_{2i+3}}{\d_{2i+2}}}\sm B_{\sqrt\fr{\d_{2i+1}}{\d_{2i+2}}}}|y|^{-(\b_{2i+1}+2)p}\l(1+|y|^{\b_{2i+2}}\r)^{ap}\mrm dy\\
&\le&C\l(\fr{\d_{2i+2}}{\d_{2i+3}}\r)^{\b_{2i+3}p}\d_{2i+2}^{2-2p}\l(\l(\fr{\d_{2i+1}}{\d_{2i+2}}\r)^\fr{2-(\b_{2i+1}+2)p}2+\l(\fr{\d_{2i+3}}{\d_{2i+2}}\r)^\fr{2+(\b_{2i+3}-2)p}2\r)\\
&=&C\l(\fr{\d_{2i+2}}{\d_{2i+3}}\r)^{\b_{2i+3}p}\d_{2i+2}^{2-2p}\l(\l(\fr{\d_{2i+2}}{\d_{2i+3}}\r)^\fr{(2-(\b_{2i+1}+2)p)\b_{2i+3}}{2\b_{2i+1}}+\l(\fr{\d_{2i+2}}{\d_{2i+3}}\r)^\fr{-2-(\b_{2i+3}-2)p}2\r)\\
&\le&C\d_{2i+2}^{2-2p}\l(\fr{\d_{2i+2}}{\d_{2i+3}}\r)^\fr{\b_{2i+3}(\b_{2i+1}p+2-2p)}{2\b_{2i+1}}\\
&\le&C\d_{2i+2}^{2-2p}|\la|^{\g'}\\
&\le&C|\la|^{\g_1}.
\eeqy
This argument has to be slightly modified when $k=2l+2$; in this case, none of the two products in the second line appear and therefore we have $\la_1^p$ in place of $\d_{2i+3}^{-\b_{2i+3}p}$:
\beqy
I'''_k&\le&C\la_1^p\d_k^{2+(\b_{k+1}-2)p}\int_{B_\fr{\mrm{diam}\O}{\d_k}\sm B_{\sqrt\fr{\d_{k-1}}{\d_k}}}|y|^{-(\b_{k-1}+2)p}\l(1+|y|^{\b_k}\r)^{ap}\mrm dy\\
&\le&C\la_1^p\d_k^{2+(\b_{k+1}-2)p}\l(\l(\fr{\d_{k-1}}{\d_k}\r)^\fr{2-(\b_{k-1}+2)p}2+\d_k^{-2-(\b_{k+1}-2)p}\r)\\
&\le&C\la_1^p\l(\d_k^{2+(\b_{k+1}-2)p}\l(\fr{\d_{k-1}}{\d_k}\r)^\fr{2-(\b_{k-1}+2)p}2+1\r)\\
&\le&C|\la|^{\g_1}.
\eeqy
The same argument works for $I'''_{2i+1}$, with a slight modification needed now for $I'''_{1}$: this time in the second line we do not have any of the sums in the power of $|x|$ and we get:
\beqy
I'''_{1}&\le&C\d_2^{-\b_2p}\int_{A_1}|x|^{(\a_2-2)p}\l(\d_1^{\b_1}+|x|^{\b_1}\r)^{bp}\mrm dx\\
&\le&C\l(\fr{\d_1}{\d_2}\r)^{\b_2p}\d_2^{2-2p}\int_{B_{\sqrt\fr{\d_2}{\d_1}}}|y|^{(\a_2-2)p}\l(1+|y|^{\b_1}\r)^{bp}\mrm dy\\
&\le&C\l(\fr{\d_1}{\d_2}\r)^\fr{(\b_2+2)p-2}2\d_2^{2-2p}\\
&\le&C|\la|^{\g_1};\eeqy
this concludes the proof.
\epf\

\subsection{The linear operator $\mcal S_\la$}

\blem\label{sla}${}$\\
Let $\mcal S_\la$ be defined by \eqref{slambda}.\\
There exists $p_0>1$ and $\g_2>0$ such that for any $p\in[1,p_0)$
$$\|\mcal S_\la\phi\|_p=O(|\la|^{\g_2}\|\phi\|).$$
\elem\

\bpf${}$\\
We can estimate $\|\mcal S_\la\|_p$ by arguing as in Lemma \ref{rla}:
\beqa
\nonumber&&\int_\O\l|\sum_{j=0}^{\l[\fr{k-1}2\r]}|\cd|^{\b_{2j+1}}e^{w_{2j+1}}-2\la_1h_1e^{W_{1,\la}}\r|^p\\
\nonumber&=&\int_\O\l|\sum_{j=0}^{\l[\fr{k-1}2\r]}|\cd|^{\b_{2j+1}-2}e^{w_{2j+1}}-2\la_1h_1e^{\sum_{m=0}^{\l[\fr{k-1}2\r]}\mrm Pw_{2m+1}-\fr{a}2\sum_{m=0}^{\l[\fr{k-2}2\r]}\mrm Pw_{2j+2}}\r|^p\\
\nonumber&\le&C\sum_{j=0}^{\l[\fr{k-1}2\r]}\int_{\mcal A_{2j+1}}\l||\cd|^{\b_{2j+1}-2}e^{w_{2j+1}}-2\la_1h_1e^{\sum_{m=0}^{\l[\fr{k-1}2\r]}\mrm Pw_{2m+1}-\fr{a}2\sum_{m=0}^{\l[\fr{k-2}2\r]}\mrm Pw_{2m+2}}\r|^p\\
\nonumber&+&C\sum_{j=0}^{\l[\fr{k-1}2\r]}\sum_{i=1,i\ne 2j+1}^k\int_{\mcal A_i}\l||\cd|^{\b_{2j+1}-2}e^{w_{2j+1}}\r|^p+C\sum_{i=0}^{\l[\fr{k-2}2\r]}\int_{\mcal A_{2i+2}}\l|\la_1h_1e^{\sum_{m=0}^{\l[\fr{k-1}2\r]}\mrm Pw_{2m+1}-\fr{a}2\sum_{m=0}^{\l[\fr{k-2}2\r]}\mrm Pw_{2m+2}}\r|^p\\
\nonumber&\le&C\l(\sum_{j=0}^{\l[\fr{k-1}2\r]}I'_{2j+1}+\sum_{j=0}^{\l[\fr{k-1}2\r]}\sum_{i=1,i\ne 2j+1}^kI''_{i,2j+1}+\sum_{i=0}^{\l[\fr{k-2}2\r]}I'''_{2i+2}\r)\\
\label{w}&\le&C|\la|^{\g_1};
\eeqa
and the same estimates also work the other components of $\mcal S_\la$.\\
Then we suffice to apply H\"older and Sobolev inequalities, with $q$ so close to $1$ that the previous estimates hold for $\|S_\la\|_{pq}$:
$$\|\mcal S_\la\phi\|_p\le\|\mcal S_\la\|_{pq}\|\phi\|_\fr{pq}{q-1}\le C|\la|^{\g_2}\|\phi\|$$
\epf\

\subsection{The quadratic term $\mcal N_\la$}

\blem\label{nla}${}$\\
Let $\mcal N_\la$ be defined by \eqref{nlambda}.\\
There exists $p_0>1,C>0$ and $\g_3>0$ such that for any $p\in[1,p_0)$
\bequ\label{n1}
\|\mcal N_\la(\phi)-\mcal N_\la(\phi')\|_p=O\l(|\la|^{\g_3(1-p)}\|\phi-\phi'\|(\|\phi\|+\|\phi'\|)e^{C\l(\|\phi\|^2+\l\|\phi'\r\|^2\r)}\r),
\eequ
and in particular
\bequ\label{n2}
\|\mcal N_\la(\phi)\|_p=O\l(|\la|^{\g_3(1-p)}\|\phi\|^2e^{C\|\phi\|^2}\r).
\eequ
\elem\

\bpf${}$\\
By writing
$$\mcal N_\la(\phi)-\mcal N_\la(\phi')=\l(\barr{cc}2&-a\\-b&2\earr\r)\l(\barr{cc}\la_1h_1e^{W_{1,\la}}\l(e^{\phi_1}-e^{\phi_1'}-\phi_1+\phi_1'\r)\\\la_2h_2e^{W_{2,\la}}\l(e^{\phi_2}-e^{\phi_2'}-\phi_2+\phi_2'\r)\earr\r),$$
we suffice to provide $L^p$ estimates for $\la_ih_ie^{W_{i,\la}}\l(e^{\phi_i}-e^{\phi'_i}-\phi_i+\phi'_i\r)$ for $i=1,2$.\\
By the elementary inequality
$$\l|e^t-e^s-t+s\r|\le|t-s|(|t|+|s|)e^{|t|+|s|}\q\q\q\fa\,t,s\in\R$$
and H\"older, Sobolev and Moser-Trudinger inequalities we get
\beqy
&&\int_\O\l|\la_ih_ie^{W_{i,\la}}\l(e^{\phi_i}-e^{\phi'_i}-\phi_i+\phi'_i\r)\r|^p\\
&\le&\int_\O\l|\la_ih_ie^{W_{i,\la}}\r|^p|\phi_i-\phi'_i|^p\l(|\phi_i|^p+\l|\phi'_i\r|^p\r)e^{p\l(|\phi_i|+\l|\phi'_i\r|\r)}\\
&\le&\l(\int_\O\l|\la_ih_ie^{W_{i,\la}}\r|^{pq}\r)^\fr{1}q\l(\int_\O|\phi_i-\phi'_i|^{ps}\r)^\fr{1}s\l(\l(\int_\O|\phi_i|^{ps}\r)^\fr{1}s+\l(\int_\O|\phi'_i|^{ps}\r)^\fr{1}s\r)\l(\int_\O e^{\fr{pqrs}{qrs-qr-qs-rs}\l(|\phi_i|+\l|\phi'_i\r|\r)}\r)^{1-\fr{1}q-\fr{1}r-\fr{1}s}\\
&\le&C\l(\int_\O\l|\la_ih_ie^{W_{i,\la}}\r|^{pq}\r)^\fr{1}q\|\phi_i-\phi'_i\|\l(\|\phi_i\|+\l\|\phi'_i\r\|\r)e^{\fr{p^2qrs}{qrs-qr-qs-rs}\l(\|\phi_i\|+\|\phi'_i\|\r)^2};
\eeqy
therefore, we just have to estimate $\la_ih_ie^{W_{i,\la}}$ in $L^p(\O)$.\\
The computations from Lemma \ref{rla} and \eqref{suptheta} yield:
\beqa
\nonumber\int_\O\l|\la_ih_ie^{W_{i,\la}}\r|^p&\le&C\sum_{j=0}^{\l[\fr{k-i}2\r]}\int_{\mcal A_{2j+i}}\l||x|^{\b_{2j+i}-2}e^{w_{2j+i}(x)+\Th_{2j+i}\l(\fr{x}{\d_{2j+i}}\r)}\r|^p\mrm dx+C\sum_{j=0}^{\l[\fr{k-3+i}2\r]}I'''_{2j-3+i}\\
\nonumber&\le&\b_{2j+i}^p\sum_{j=0}^{\l[\fr{k-i}2\r]}\d_{2j+i}^{2-2p}\int_{\fr{\mcal A_{2j+i}}{\d_{2j+i}}}\fr{|y|^{(\b_{2j+i}-2)p}}{(1+|y|^{\b_{2j+i}})^{2p}}e^{p|\Th_{2j+i}(y)|}\mrm dy+o(1)\\
\nonumber&\le&C\sum_{j=0}^{\l[\fr{k-i}2\r]}\d_{2j+i}^{2-2p}\int_{\fr{\mcal A_{2j+i}}{\d_{2j+i}}}\fr{|y|^{(\b_{2j+i}-2)p}}{(1+|y|^{\b_{2j+i}})^{2p}}\mrm dy+o(1)\\
\nonumber&\le&C\sum_{j=0}^{\l[\fr{k-i}2\r]}\d_{2j+i}^{2-2p}\\
\label{lp}&\le&C|\la|^{\g_3(1-p)}
\eeqa
hence \eqref{n1} is proved.\\
\eqref{n2} just follows from \eqref{n1} after setting $\phi'=0$.
\epf\

\section{Linear theory}\label{linear}\

In this section we develop a linear theory for the linear operator $\mcal L_\la$ defined in \eqref{llambda}.\\
The following proposition, whose proof will take up the whole section, is inspired by \cite{gp} (Proposition $4.1$) and \cite{mpw} (Proposition $4.1$).\\

\bprop\label{lla}${}$\\
For any $p>1$ there exists $\ol\la>0$ and $C>0$ such that for any $\la\in\l(0,\ol\la\r)\x\l(0,\ol\la\r)$ and any $\psi\in\mbf H
$ there exists a unique $\phi\in\mbf H$ solution of
$$\mcal L_\la\phi=\psi\q\tx{on }\O,$$
satisfying
$$\|\phi\|\le C\log\fr{1}{|\la|}\|\psi\|_p$$
\eprop\

\bpf${}$\\
Suppose the statement is not true. This means that there exist $p>1$ and sequences $\{\la_n\}_{n\in\N}\sub\R_{>0}^2,\,\{\psi_n\}_{n\in\N}\sub\mbf H,\,\{\phi_n\}_{n\in\N}\sub\mbf H$ such that
\bequ\label{phn}
\l\{\bl-\D\phi_{n,1}-\sum_{j=0}^{\l[\fr{k-1}2\r]}2\b_{2j+1}^2\fr{\d_{n,2j+1}^{\b_{2j+1}}|\cd|^{\b_{2j+1}-2}}{\l(\d_{n,2j+1}^{\b_{2j+1}}+|\cd|^{\b_{2j+1}}\r)^2}\phi_{n,1}+\fr{a}2\sum_{j=0}^{\l[\fr{k-2}2\r]}2\b_{2j+2}^2\fr{\d_{n,2j+2}^{\b_{2j+2}}|\cd|^{\b_{2j+2}-2}}{\l(\d_{n,2j+2}^{\b_{2j+2}}+|\cd|^{\b_{2j+2}}\r)^2}\phi_{n,2}=\psi_{n,1}\\-\D\phi_{n,2}-\sum_{j=0}^{\l[\fr{k-2}2\r]}2\b_{2j+1}^2\fr{\d_{n,2j+2}^{\b_{2j+2}}|\cd|^{\b_{2j+2}-2}}{\l(\d_{n,2j+2}^{\b_{2j+2}}+|\cd|^{\b_{2j+2}}\r)^2}\phi_{n,2}+\fr{b}2\sum_{j=0}^{\l[\fr{k-1}2\r]}2\b_{2j+1}^2\fr{\d_{n,2j+1}^{\b_{2j+1}}|\cd|^{\b_{2j+1}-2}}{\l(\d_{n,2j+1}^{\b_{2j+1}}+|\cd|^{\b_{2j+1}}\r)^2}\phi_{n,1}=\psi_{n,2}\\\la_n\us{n\to+\infty}\to0\q\q\q\q\q\q\q\q\q\q\|\phi_n\|=1\q\q\q\q\q\q\q\q\q\q\log\fr{1}{|\la_n|}\|\psi_n\|_p\us{n\to+\infty}\to0\earr\r.,
\eequ
where $\d_{n,\ell}$ is defined as in \eqref{deltaj} with $\la_{n,1},\la_{n,2}$ in place of $\la_1,\la_2$.\\
We will divide the proof in six steps.\\

\bite
\item[\emph{Step $1$:}]\emph{$\int_\O\fr{\d_{n,2j+i}^{\b_{2j+i}}|\cd|^{\b_{2j+i}-2}}{\l(\d_{n,2j+i}^{\b_{2j+i}}+|\cd|^{\b_{2j+i}}\r)^2}|\phi_{n,i}|^2=O(1)$ for any $i=1,2$, $j=0,\ds,\l[\fr{k-i}2\r]$.}\\
If we multiply both sides of the first equation in \eqref{phn} by $\phi_{n,1}$ and both sides of the second equation by $\fr{a}2\phi_{n,1}$ and then we sum the two equalities we get
\beqy
&&\l(1-\fr{ab}4\r)\sum_{j=0}^{\l[\fr{k-1}2\r]}2\b_{2j+1}^2\int_\O\fr{\d_{n,2j+1}^{\b_{2j+1}}|\cd|^{\b_{2j+1}-2}}{\l(\d_{n,2j+1}^{\b_{2j+1}}+|\cd|^{\b_{2j+1}}\r)^2}|\phi_{n,1}|^2\\
&=&\int_\O|\n\phi_{n,1}|^2-\int_\O\psi_{n,1}\phi_{n,1}+\fr{a}2\int_\O\n\phi_{n,1}\cd\n\phi_{n,2}-\fr{a}2\int_\O\psi_{n,2}\phi_{n,1}\\
&\le &C\l(\|\phi_n\|^2+\|\psi_n\|_p\|\phi_n\|\r)\\
&\le&C;
\eeqy
similarly, by multiplying the first equation in \eqref{phn} by $\fr{b}2\phi_{n,1}$, the second equation by $\phi_{n,2}$ and then summing, we get
$$\l(1-\fr{ab}4\r)\sum_{j=0}^{\l[\fr{k-2}2\r]}2\b_{2j+2}^2\int_\O\fr{\d_{n,2j+2}^{\b_{2j+2}}|\cd|^{\b_{2j+2}-2}}{\l(\d_{n,2j+2}^{\b_{2j+2}}+|\cd|^{\b_{2j+2}}\r)^2}|\phi_{n,2}|^2=O(1).$$
Therefore, the claim is proved if $ab\ne 4$.\\
On the other hand, if $ab=4$, then summing the first equation in \eqref{phn} and the second equation multiplied by $\fr{a}2=\fr{2}b$ gives
$$\l\{\bll-\D\l(\phi_{n,1}+\fr{a}2\phi_{n,2}\r)=\psi_{n,1}+\fr{a}2\psi_{n,2}&\tx{in }\O\\\phi_{n,1}+\fr{a}2\phi_{n,2}=0&\tx{on }\pa\O\earr\r.,$$
hence standard regularity theory yields
$$\l\|\phi_{n,1}+\fr{a}2\phi_{n,2}\r\|_\infty\le C\l\|\psi_{n,1}+\fr{a}2\psi_{n,2}\r\|_p=o(1).$$
Therefore, multiplying the first equation in \eqref{phn} by $\phi_{n,1}$, the second equation by $\fr{a^2}4\phi_{n,2}=\fr{4}{b^2}\phi_{n,2}$ and then summing we get
\beqa
\label{ab4}&&2\sum_{j=0}^{\l[\fr{k-1}2\r]}2\b_{2j+1}^2\int_\O\fr{\d_{n,2j+1}^{\b_{2j+1}}|\cd|^{\b_{2j+1}-2}}{\l(\d_{n,2j+1}^{\b_{2j+1}}+|\cd|^{\b_{2j+1}}\r)^2}|\phi_{n,1}|^2+\fr{a^2}2\sum_{j=0}^{\l[\fr{k-2}2\r]}2\b_{2j+2}^2\int_\O\fr{\d_{n,2j+2}^{\b_{2j+2}}|\cd|^{\b_{2j+2}-2}}{\l(\d_{n,2j+2}^{\b_{2j+2}}+|\cd|^{\b_{2j+2}}\r)^2}|\phi_{n,2}|^2\\
\nonumber&=&\int_\O\l(\sum_{j=0}^{\l[\fr{k-1}2\r]}2\b_{2j+1}^2\fr{\d_{n,2j+1}^{\b_{2j+1}}|\cd|^{\b_{2j+1}-2}}{\l(\d_{n,2j+1}^{\b_{2j+1}}+|\cd|^{\b_{2j+1}}\r)^2}\phi_{n,1}+\fr{a}2\sum_{j=0}^{\l[\fr{k-2}2\r]}2\b_{2j+2}^2\fr{\d_{n,2j+2}^{\b_{2j+2}}|\cd|^{\b_{2j+2}-2}}{\l(\d_{n,2j+2}^{\b_{2j+2}}+|\cd|^{\b_{2j+2}}\r)^2}\phi_{n,2}\r)\l(\phi_{n,1}+\fr{a}2\phi_{n,2}\r)\\
\nonumber&+&\int_\O|\n\phi_{n,1}|^2-\int_\O\psi_{n,1}\phi_{n,1}+\fr{a^2}4\int_\O|\n\phi_{n,2}|^2-\fr{a^2}4\int_\O\psi_{n,2}\phi_{n,2}\\
\nonumber&\le&\l(\sum_{j=0}^{\l[\fr{k-1}2\r]}2\b_{2j+1}^2\int_\O\fr{\d_{n,2j+1}^{\b_{2j+1}}|\cd|^{\b_{2j+1}-2}}{\l(\d_{n,2j+1}^{\b_{2j+1}}+|\cd|^{\b_{2j+1}}\r)^2}|\phi_{n,1}|+\fr{a}2\sum_{j=0}^{\l[\fr{k-2}2\r]}2\b_{2j+2}^2\int_\O\fr{\d_{n,2j+2}^{\b_{2j+2}}|\cd|^{\b_{2j+2}-2}}{\l(\d_{n,2j+2}^{\b_{2j+2}}+|\cd|^{\b_{2j+2}}\r)^2}|\phi_{n,2}|\r)^\fr{1}2\l\|\phi_{n,1}+\fr{a}2\phi_{n,2}\r\|_\infty+C\\
\nonumber&\le&o(1)\l(\sum_{j=0}^{\l[\fr{k-1}2\r]}2\b_{2j+1}^2\int_\O\fr{\d_{n,2j+1}^{\b_{2j+1}}|\cd|^{\b_{2j+1}-2}}{\l(\d_{n,2j+1}^{\b_{2j+1}}+|\cd|^{\b_{2j+1}}\r)^2}|\phi_{n,1}|^2+\fr{a}2\sum_{j=0}^{\l[\fr{k-2}2\r]}2\b_{2j+2}^2\int_\O\fr{\d_{n,2j+2}^{\b_{2j+2}}|\cd|^{\b_{2j+2}-2}}{\l(\d_{n,2j+2}^{\b_{2j+2}}+|\cd|^{\b_{2j+2}}\r)^2}|\phi_{n,2}|^2\r)^\fr{1}2+C;
\eeqa
therefore, $\int_\O\fr{\d_{n,2j+i}^{\b_{2j+i}}|\cd|^{\b_{2j+i}-2}}{\l(\d_{n,2j+i}^{\b_{2j+i}}+|\cd|^{\b_{2j+i}}\r)^2}|\phi_{n,i}|^2\le C$ for all $i,\ell$.\\

\item[\emph{Step $2$:}]\emph{The sequence $\wt\phi_{n,\ell}$, defined by $\wt\phi_{n,2j+i}(y):=\phi_{n,i}(\d_{n,2j+i}y)$, converges to $\mu_\ell\fr{1-|\cd|^{\b_\ell}}{1+|\cd|^{\b_\ell}}$ for some $\mu_\ell$, weakly in $H_{\b_\ell}\l(\R^2\r)$ and strongly in $L_{\b_\ell}\l(\R^2\r)$, where
$$\barr{rcll}L_\b\l(\R^2\r)&:=&\l\{u\in L^2_{\mrm{loc}}\l(\R^2\r):\,\fr{|\cd|^\fr{\b-2}2}{1+|\cd|^\b}u\in L^2\l(\R^2\r)\r\},&\|u\|_{L_\b}:=\l\|\fr{|\cd|^\fr{\b-2}2}{1+|\cd|^\b}u\r\|_{L^2(\R^2)};\\
H_\b\l(\R^2\r)&:=&\l\{u\in H^1_{\mrm{loc}}\l(\R^2\r):\,|\n u|+\fr{|\cd|^\fr{\b-2}2}{1+|\cd|^\b}u\in L^2\l(\R^2\r)\r\},&\|u\|_{H_\b}:=\l(\|\n u\|_{L^2(\R^2)}^2+\l\|\fr{|\cd|^\fr{\b-2}2}{1+|\cd|^\b}u\r\|_{L^2(\R^2)}^2\r)^\fr{1}2.\earr$$}
First of all, because of Step $1$, $\wt\phi_{n,\ell}$ is bounded in $H_{\b_\ell}\l(\R^2\r)$:
\beqy
\int_\fr{\O}{\d_{n,2j+i}}\l|\n\wt\phi_{n,2j+i}(y)\r|^2\mrm dy&=&\d_{n,2j+i}^2\int_\fr{\O}{\d_{n,2j+i}}|\n\phi_{n,i}(\d_{n,2j+i}y)|^2\mrm dy=\int_\O|\n\phi_{n,i}(x)|^2\mrm dx=1\\
\int_\fr{\O}{\d_{n,2j+i}}\fr{|y|^{\b_{2j+i}-2}}{(1+|y|^{\b_{2j+i}})^2}\l|\wt\phi_{n,2j+i}(y)\r|^2\mrm dy&=&\int_\O\fr{\d_{n,2j+i}^{\b_{2j+i}}|x|^{\b_{2j+i}-2}}{\l(\d_{n,2j+i}^{\b_{2j+i}}+|x|^{\b_{2j+i}}\r)^2}|\phi_{n,i}(x)|^2\mrm dx=O(1).
\eeqy
Therefore, $\wt\phi_{n,\ell}\us{n\to+\infty}\wk\wt\phi_\ell$ in $H_{\b_\ell}\l(\R^2\r)$ for some $\wt\phi_\ell\in H_{\b_\ell}\l(\R^2\r)$; moreover, the embedding $L_{\b_\ell}\l(\R^2\r)\inc H_{\b_\ell}\l(\R^2\r)$ is compact (see \cite{gp}, Proposition $6.1$; the result is stated only for $\a\ge2$ but the same argument works for any $\a>0$). From this we get $\wt\phi_{n,\ell}\us{n\to+\infty}\to\wt\phi_\ell$ in $L_{\b_\ell}\l(\R^2\r)$.\\
$\wt\phi_{n,\ell}$ solves
$$\l\{\bll-\D\wt\phi_{n,\ell}=2\b_\ell^2\fr{|\cd|^{\b_\ell-2}}{(1+|\cd|^{\b_\ell})^2}\wt\phi_{n,\ell}+\rho_{n,\ell}&\tx{in }\fr{\O}{\d_{n,\ell}}\\\wt\phi_{n,\ell}=0&\tx{on }\pa\l(\fr{\O}{\d_{n,\ell}}\r)\earr\r.,$$
where
\beqy
\rho_{n,2j+1}(y)&:=&\sum_{i=0,i\ne j}^{\l[\fr{k-1}2\r]}2\b_{2i+1}^2\fr{\d_{n,2i+1}^{\b_{2i+1}}\d_{n,2j+1}^{\b_{2i+1}}|y|^{\b_{2i+1}-2}}{\l(\d_{n,2i+1}^{\b_{2i+1}}+\d_{n,2j+1}^{\b_{2i+1}}|y|^{\b_{2i+1}}\r)^2}\phi_{n,1}(\d_{n,2j+1}y)\\
&-&\fr{a}2\sum_{i=0}^{\l[\fr{k-2}2\r]}2\b_{2i+2}^2\fr{\d_{n,2i+2}^{\b_{2i+2}}\d_{n,2j+1}^{\b_{2i+2}}|y|^{\b_{2i+2}-2}}{\l(\d_{n,2i+2}^{\b_{2i+2}}+\d_{n,2j+1}^{\b_{2i+2}}|y|^{\b_{2i+2}}\r)^2}\phi_{n,2}(\d_{n,2j+1}y)+\d_{n,2j+1}^2\psi_{n,1}(\d_{n,2j+1}y);\\
\rho_{n,2j+2}(y)&:=&\sum_{i=0,i\ne j}^{\l[\fr{k-2}2\r]}2\b_{2i+2}^2\fr{\d_{n,2i+2}^{\b_{2i+2}}\d_{n,2j+2}^{\b_{2i+2}}|y|^{\b_{2i+2}-2}}{\l(\d_{n,2i+2}^{\b_{2i+2}}+\d_{n,2j+2}^{\b_{2i+2}}|y|^{\b_{2i+2}}\r)^2}\phi_{n,2}(\d_{n,2j+2}y)\\
&-&\fr{b}2\sum_{i=0}^{\l[\fr{k-1}2\r]}2\b_{2i+1}^2\fr{\d_{n,2i+1}^{\b_{2i+1}}\d_{n,2j+2}^{\b_{2i+1}}|y|^{\b_{2i+1}-2}}{\l(\d_{n,2i+1}^{\b_{2i+1}}+\d_{n,2j+2}^{\b_{2i+1}}|y|^{\b_{2i+1}}\r)^2}\phi_{n,1}(\d_{n,2j+2}y)+\d_{n,2j+2}^2\psi_{n,2}(\d_{n,2j+2}y).
\eeqy
Let us show that $\rho_{n,\ell}\us{n\to+\infty}\to0$ in $L^1_{\mrm{loc}}\l(\R^2\sm\{0\}\r)$.\\
Any compact set $\mcal K\Sub\R^2\sm\{0\}$ will be contained, for large $n$, in $\fr{A_{n,\ell}}{\d_{n,\ell}}:=\l\{y\in\fr{\O}{\d_{n,\ell}}:\,\sqrt{\fr{\d_{n,\ell-1}}{\d_{n,\ell}}}\le|y|\le\sqrt{\fr{\d_{n,\ell}}{\d_{n,\ell+1}}}\r\}$; therefore, by the estimate \eqref{i2},
\beqy
\int_{\mcal K}|\rho_{n,\ell}|&\le&\int_\fr{A_{n,\ell}}{\d_{n,\ell}}|\rho_{n,\ell}(y)|\mrm dy\\
&\le&C\sum_{i=0,i\ne j}^k\int_\fr{A_{n,\ell}}{\d_{n,\ell}}\fr{\d_{n,i}^{\b_i}\d_{n,\ell}^{\b_i}|y|^{\b_i-2}}{\l(\d_{n,i}^{\b_i}+\d_{n,\ell}^{\b_i}|y|^{\b_i}\r)^2}(|\phi_{n,1}(\d_{n,\ell}y)|+|\phi_{n,2}(\d_{n,\ell}y)|)\mrm dy\\
&+&\d_{n,\ell}^2\int_\fr{\O}{\d_{n,\ell}}(|\psi_{n,1}(\d_{n,\ell}y)|+|\psi_{n,2}(\d_{n,\ell}y)|)\mrm dy\\
&=&C\sum_{i=0,i\ne j}^k\int_{A_{n,\ell}}\fr{\d_{n,i}^{\b_i}|x|^{\b_i-2}}{\l(\d_{n,i}^{\b_i}+|x|^{\b_i}\r)^2}(|\phi_{n,1}(x)|+|\phi_{n,2}(x)|)\mrm dx+C\int_\O(|\psi_{n,1}(x)|+|\psi_{n,2}(x)|)\mrm dx\\
&\le&C\l(\sum_{i=0,i\ne j}^k\int_{A_{n,\ell}}\l|\fr{\d_{n,i}^{\b_i}|x|^{\b_i-2}}{\l(\d_{n,i}^{\b_i}+|x|^{\b_i}\r)^2}\r|^q\mrm dx\r)^\fr{1}q\|\phi_n\|_\fr{q}{q-1}+C\|\psi_n\|_p\\
&\le&C|\la_n|^{\g_1}\|\phi_n\|+C\|\psi_n\|_p\\
&\us{n\to+\infty}\to&0.
\eeqy
Therefore, the weak limit $\wt\phi_\ell$ must be a solution of
$$-\D\wt\phi_\ell=2\b_\ell^2\fr{|\cd|^{\b_\ell-2}}{(1+|\cd|^{\b_\ell})^2}\wt\phi_\ell\q\tx{in }\R^2\sm\{0\}.$$
Finally, by the properties of weak convergence we get $\int_{\R^2}\l|\n\wt\phi_\ell\r|^2\le1$, therefore $\wt\phi_\ell$ must be a solution on the whole plane; by Proposition \ref{intere} we get $\wt\phi_\ell=\mu_\ell\fr{1-|\cd|^{\b_\ell}}{1+|\cd|^{\b_\ell}}$.\\

\item[\emph{Step $3$:}]\emph{$\s_{n,\ell}:=\log\fr{1}{|\la_n|}\int_\fr{\O}{\d_{n,\ell}}2\b_\ell^2\fr{|\cd|^{\b_\ell-2}}{\l(1+|\cd|^{\b_\ell}\r)^2}\wt\phi_{n,\ell}\us{n\to+\infty}\to0$ for all $\ell$'s.}\\
Define $Z_{n,\ell}:=\fr{\d_{n,\ell}^{\b_\ell}-|\cd|^{\b_\ell}}{\d_{n,\ell}^{\b_\ell}+|\cd|^{\b_\ell}}$, which solves (see Theorem \ref{intere})
$$-\D Z_{n,\ell}=2\b_\ell^2\fr{\d_{n,\ell}^{\b_\ell}|\cd|^{\b_\ell-2}}{\l(\d_{n,\ell}^{\b_\ell}+|\cd|^{\b_\ell}\r)^2}Z_{n,\ell}\q\tx{in }\R^2;$$
consider now its projection $\mrm PZ_{n,\ell}$ on $H^1_0(\O)$, namely (see \eqref{pro}) the solution of
\bequ\label{eqpz}
\l\{\bll-\D(\mrm PZ_{n,\ell})=2\b_\ell^2\fr{\d_{n,\ell}^{\b_\ell}|\cd|^{\b_\ell-2}}{\l(\d_{n,\ell}^{\b_\ell}+|\cd|^{\b_\ell}\r)^2}Z_{n,\ell}&\tx{in }\O\\\mrm PZ_{n,\ell}=0&\tx{on }\pa\O\earr\r.
\eequ
As in Lemma \ref{pw}, the maximum principle gives
\bequ\label{pz}
\mrm PZ_{n,\ell}=Z_{n,\ell}+1+O\l(\d_{n,\ell}^{\b_\ell}\r)=\fr{2\d_{n,\ell}^{\b_\ell}}{\d_{n,\ell}^{\b_\ell}+|\cd|^{\b_\ell}}+O\l(\d_{n,\ell}^{\b_\ell}\r),
\eequ
hence
\bequ\label{pzd}
\mrm PZ_{n,i}(\d_{n,\ell}y)=\l\{\bll2\fr{\l(\fr{\d_{n,i}}{\d_{n,\ell}}\r)^{\b_i}}{\l(\l(\fr{\d_{n,i}}{\d_{n,\ell}}\r)^{\b_i}+|y|^{\b_i}\r)}+O\l(\d_{n,i}^{\b_i}\r)&\tx{if }i<\ell\\\fr{2}{1+|y|^{\b_i}}+O\l(\d_{n,i}^{\b_i}\r)&\tx{if }i=\ell\\2-2\fr{\l(\fr{\d_{n,\ell}}{\d_{n,i}}\r)^{\b_i}|y|^{\b_i}}{\l(1+\l(\fr{\d_{n,\ell}}{\d_{n,i}}\r)^{\b_i}|y|^{\b_i}\r)}+O\l(\d_{n,i}^{\b_i}\r)&\tx{if }i>\ell\earr\r.
\eequ
Recall now the first equation of \eqref{phn} and multiply it by $\log\fr{1}{|\la_n|}\mrm PZ_{n,2i+1}$; then, multiply by $\log\fr{1}{|\la_n|}\phi_{n,1}$ the equation \eqref{eqpz} satisfied by $\mrm PZ_{n,2i+1}$ and subtract the two quantities: we get
\beqy
0&=&\ub{\log\fr{1}{|\la_n|}\int_\O2\b_{2i+1}^2\fr{\d_{n,2i+1}^{\b_{2i+1}}|\cd|^{\b_{2i+1}-2}}{\l(\d_{n,2i+1}^{\b_{2i+1}}+|\cd|^{\b_{2i+1}}\r)^2}\phi_{n,1}(\mrm PZ_{n,2i+1}-Z_{n,2i+1})}_{=:I'_{n,2i+1}}\\
&+&\sum_{j=0,j\ne i}^{\l[\fr{k-1}2\r]}\ub{\log\fr{1}{|\la_n|}\int_\O2\b_{2j+1}^2\fr{\d_{n,2j+1}^{\b_{2j+1}}|\cd|^{\b_{2j+1}-2}}{\l(\d_{n,2j+1}^{\b_{2j+1}}+|\cd|^{\b_{2j+1}}\r)^2}\phi_{n,1}\mrm PZ_{n,2i+1}}_{=:I''_{n,2i+1,2j+1}}\\
&-&\fr{a}2\sum_{j=0}^{\l[\fr{k-2}2\r]}\ub{\log\fr{1}{|\la_n|}\int_\O2\b_{2j+2}^2\fr{\d_{n,2j+2}^{\b_{2j+2}}|\cd|^{\b_{2j+2}-2}}{\l(\d_{n,2j+2}^{\b_{2j+2}}+|\cd|^{\b_{2j+2}}\r)^2}\phi_{n,2}\mrm PZ_{n,2i+1}}_{=:I''_{n,2i+1,2j+2}}+\ub{\log\fr{1}{|\la_n|}\int_\O\psi_{n,1}\mrm PZ_{n,2i+1}}_{=:I'''_{n,2i+1}}.
\eeqy
To estimate $I'_{n,2i+1}$ we use \eqref{pz}, then the boundedness in $L_{\b_{2i+1}}\l(\R^2\r)$ and the definitions of $\d_{n,i}$:
\beqy
I'_{n,2i+1}&=&\log\fr{1}{|\la_n|}\int_\fr{\O}{\d_{n,2i+1}}2\b_{2i+1}^2\fr{|y|^{\b_{2i+1}-2}}{(1+|y|^{\b_{2i+1}})^2}\wt\phi_{n,2i+1}(y)(\mrm PZ_{n,2i+1}(\d_{n,2i+1}y)-Z_{n,2i+1}(\d_{n,2i+1}y))\mrm dy\\
&=&\log\fr{1}{|\la_n|}\int_\fr{\O}{\d_{n,2i+1}}2\b_{2i+1}^2\fr{|y|^{\b_{2i+1}-2}}{(1+|y|^{\b_{2i+1}})^2}\wt\phi_{n,2i+1}(y)\mrm dy\\
&+&O\l(\d_{n,2i+1}^{\b_{2i+1}}\log\fr{1}{|\la_n|}\int_\fr{\O}{\d_{n,2i+1}}2\b_{2i+1}^2\fr{|y|^{\b_{2i+1}-2}}{(1+|y|^{\b_{2i+1}})^2}\l|\wt\phi_{n,2i+1}(y)\r|\mrm dy\r)\\
&=&\s_{n,2i+1}+O\l(\d_{n,2i+1}^{\b_{2i+1}}\log\fr{1}{|\la_n|}\l\|\wt\phi_{n,2i+1}\r\|_{L_{\b_{2i+1}}}\r)\\
&=&\s_{n,2i+1}+o(1).
\eeqy
Concerning the terms in $I''_{n,2i+1,2j+1}$, we proceed differently depending whether $j<i$ or $j>i$: in the former case, using \eqref{pzd} and choosing $q$ very close to $1$ we get
\beqy
I''_{n,2i+1,2j+1}&=&\log\fr{1}{|\la_n|}\int_\fr{\O}{\d_{n,2j+1}}2\b_{2j+1}^2\fr{|y|^{\b_{2j+1}-2}}{(1+|y|^{\b_{2j+1}})^2}\wt\phi_{n,2j+1}(y)\mrm PZ_{n,2i+1}(\d_{n,2j+1}y)\mrm dy\\
&=&2\log\fr{1}{|\la_n|}\int_\fr{\O}{\d_{n,2j+1}}2\b_{2j+1}^2\fr{|y|^{\b_{2j+1}-2}}{(1+|y|^{\b_{2j+1}})^2}\wt\phi_{n,2j+1}(y)\mrm dy\\
&-&2\l(\fr{\d_{n,2j+1}}{\d_{n,2i+1}}\r)^{\b_{2i+1}}\log\fr{1}{|\la_n|}\int_\fr{\O}{\d_{n,2j+1}}2\b_{2j+1}^2\fr{|y|^{\b_{2j+1}+\b_{2i+1}-2}}{(1+|y|^{\b_{2j+1}})^2\l(1+\l(\fr{\d_{n,2j+1}}{\d_{n,2i+1}}\r)^{\b_{2i+1}}|y|^{\b_{2i+1}}\r)}\wt\phi_{n,2j+1}(y)\mrm dy\\
&+&O\l(\d_{n,2i+1}^{\b_{2i+1}}\log\fr{1}{|\la_n|}\int_\fr{\O}{\d_{n,2j+1}}2\b_{2j+1}^2\fr{|y|^{\b_{2j+1}-2}}{(1+|y|^{\b_{2j+1}})^2}\l|\wt\phi_{n,2j+1}(y)\r|\mrm dy\r)\\
&=&2\s_{n,2j+1}+o(1)\\
&+&O\l(\l(\fr{\d_{n,2j+1}}{\d_{n,2i+1}}\r)^{\b_{2i+1}}\log\fr{1}{|\la_n|}\l(\int_{\R^2}\l|\fr{|y|^{\b_{2j+1}+\b_{2i+1}-2}}{(1+|y|^{\b_{2j+1}})^2\l(1+\l(\fr{\d_{n,2j+1}}{\d_{n,2i+1}}\r)^{\b_{2i+1}}|y|^{\b_{2i+1}}\r)}\r|^q\mrm dy\r)^\fr{1}q\l\|\wt\phi_{n,2j+1}\r\|_\fr{q}{q-1}\r)+o(1)\\
&=&2\s_{n,2j+1}+O\l(\l(\fr{\d_{n,2j+1}}{\d_{n,2i+1}}\r)^{\b_{2i+1}}\log\fr{1}{|\la_n|}\l(\int_{\R^2\sm B_1(0)}\fr{|y|^{(\b_{2i+1}-\b_{2j+1}-2)q}}{\l(1+\l(\fr{\d_{n,2j+1}}{\d_{n,2i+1}}\r)^{\b_{2i+1}}|y|^{\b_{2i+1}}\r)^q}\r)^\fr{1}q\d_{n,2j+1}^{-2\l(1-\fr{1}q\r)}\|\phi_{n,1}\|_\fr{q}{q-1}\r)+o(1)\\
&=&2\s_{n,2j+1}+O\l(\l(\fr{\d_{n,2j+1}}{\d_{n,2i+1}}\r)^{\b_{2i+1}}\log\fr{1}{|\la_n|}\l(\fr{\d_{n,2j+1}}{\d_{n,2i+1}}\r)^{\min\l\{0,\b_{2j+1}-\b_{2i+1}+2\l(1-\fr{1}q\r)\r\}}\d_{n,2j+1}^{-2\l(1-\fr{1}q\r)}\|\phi_n\|\r)+o(1)\\
&=&2\s_{n,2j+1}+o(1);
\eeqy
in the latter case,
\beqy
I''_{n,2i+1,2j+1}&=&\log\fr{1}{|\la_n|}\int_\fr{\O}{\d_{n,2j+1}}2\b_{2j+1}^2\fr{|y|^{\b_{2j+1}-2}}{(1+|y|^{\b_{2j+1}})^2}\wt\phi_{n,2j+1}(y)\mrm PZ_{n,2i+1}(\d_{n,2j+1}y)\mrm dy\\
&=&\l(\fr{\d_{n,2i+1}}{\d_{n,2j+1}}\r)^{\b_{2i+1}}\log\fr{1}{|\la_n|}\int_\fr{\O}{\d_{n,2j+1}}2\b_{2j+1}^2\fr{|y|^{\b_{2j+1}-2}}{(1+|y|^{\b_{2j+1}})^2\l(\l(\fr{\d_{n,2i+1}}{\d_{n,2j+1}}\r)^{\b_{2i+1}}+|y|^{\b_{2i+1}}\r)}\wt\phi_{n,2j+1}(y)\mrm dy\\\\
&+&O\l(\d_{n,2i+1}^{\b_{2i+1}}\log\fr{1}{|\la_n|}\int_\fr{\O}{\d_{n,2j+1}}2\b_{2j+1}^2\fr{|y|^{\b_{2j+1}-2}}{(1+|y|^{\b_{2j+1}})^2}\l|\wt\phi_{n,2j+1}(y)\r|\mrm dy\r)\\
&=&O\l(\l(\fr{\d_{n,2i+1}}{\d_{n,2j+1}}\r)^{\b_{2i+1}}\log\fr{1}{|\la_n|}\l(\int_{B_1(0)}\fr{|y|^{(\b_{2j+1}-2)q}}{\l(\l(\fr{\d_{n,2i+1}}{\d_{n,2j+1}}\r)^{\b_{2i+1}}+|y|^{\b_{2i+1}}\r)^q}\mrm dy\r)^\fr{1}q\d_{n,2j+1}^{-2\l(1-\fr{1}q\r)}\|\phi_{n,1}\|_\fr{q}{q-1}\r)+o(1)\\
&=&O\l(\l(\fr{\d_{n,2i+1}}{\d_{n,2j+1}}\r)^{\b_{2i+1}}\log\fr{1}{|\la_n|}\l(\fr{\d_{n,2i+1}}{\d_{n,2j+1}}\r)^{\min\l\{0,\b_{2j+1}-\b_{2i+1}+2\l(1-\fr{1}q\r)\r\}}\d_{n,2j+1}^{-2\l(1-\fr{1}q\r)}\|\phi_n\|\r)+o(1)\\
&\us{n\to+\infty}\to&0.
\eeqy
The same argument shows that
$$I'''_{n,2i+1,2j+2}=\l\{\bll2\s_{n,2j+2}+o(1)&\tx{if }j<i\\o(1)&\tx{if }j\ge i\earr\r..$$
Finally, since $\|\mrm PZ_{n,\ell}\|_\infty\le C$,
$$|I'''_{n,2i+1}|\le\log\fr{1}{|\la_n|}\|\psi_{n,1}\|_1\|\mrm PZ_{n,2i+1}\|_\infty\le C\log\fr{1}{|\la_n|}\|\psi_n\|_p\us{n\to+\infty}\to0.$$
Therefore, we get:
\bequ\label{s1}
\s_{n,2i+1}+2\sum_{j=0}^{i-1}\s_{n,2j+1}-a\sum_{j=0}^{i-1}\s_{n,2j+2}=o(1);
\eequ
a similar argument yields
\bequ\label{s2}
\s_{n,2i+2}+2\sum_{j=0}^{i-1}\s_{n,2j+2}-b\sum_{j=0}^i\s_{n,2j+1}=o(1).
\eequ
Putting \eqref{s1} and \eqref{s2} together we get $\s_{n,i}=o(1)$ for all $i$'s.\\

\item[\emph{Step $4$:}]\emph{$\mu_\ell=0$ for all $j$'s.}\\
We recall the solution $\mrm Pw_{n,\ell}=\mrm Pw_{\d_{n,\ell}}^{\b_\ell}$ of
\bequ\label{pwn}
\l\{\bll-\D(\mrm Pw_{n,\ell})=2\b_\ell^2\fr{\d_{n,\ell}^{\b_\ell}|\cd|^{\b_\ell-2}}{\l(\d_{n,\ell}^{\b_\ell}+|\cd|^{\b_\ell}\r)^2}&\tx{in }\O\\\mrm Pw_{n,\ell}=0&\tx{on }\pa\O\earr\r..
\eequ
We multiply by $\mrm Pw_{n,2i+1}$ the first equation of \eqref{phn}, then we multiply by $\phi_{n,1}$ the equation \eqref{pwn} satisfied by $\mrm Pw_{n,2i+1}$; their difference gives
\beqy
0&=&\sum_{j=0}^{\l[\fr{k-1}2\r]}\ub{\int_\O2\b_{2j+1}^2\fr{\d_{n,2j+1}^{\b_{2j+1}}|\cd|^{\b_{2j+1}-2}}{\l(\d_{n,2j+1}^{\b_{2j+1}}+|\cd|^{\b_{2j+1}}\r)^2}\phi_{n,1}\mrm Pw_{n,2i+1}}_{=:J'_{n,2i+1,2j+1}}-\fr{a}2\sum_{j=0}^{\l[\fr{k-2}2\r]}\ub{\int_\O2\b_{2j+2}^2\fr{\d_{n,2j+2}^{\b_{2j+2}}|\cd|^{\b_{2j+2}-2}}{\l(\d_{n,2j+2}^{\b_{2j+2}}+|\cd|^{\b_{2j+2}}\r)^2}\phi_{n,2}\mrm Pw_{n,2i+1}}_{=:J'_{n,2i+1,2j+2}}\\
&+&\ub{\int_\O\psi_{n,1}\mrm Pw_{n,2i+1}}_{=:J''_{n,2i+1}}-\ub{\int_\O2\b_{2i+1}^2\fr{\d_{n,2i+1}^{\b_{2i+1}}|\cd|^{\b_{2i+1}-2}}{\l(\d_{n,2i+1}^{\b_{2i+1}}+|\cd|^{\b_{2i+1}}\r)^2}\phi_{n,1}}_{=:J'''_{n,2i+1}}.
\eeqy
We start by estimating $J'_{n,i,\ell}$, considering as before only the case of odd indexes.\\
For $\ell<i$ we use \eqref{pw1}, the definition of $\d_{n,\ell}$ and the vanishing of $\s_{n,\ell}$. Notice that, to handle with $\mrm Pw_i$, \eqref{pw2} would not suffice hence we need sharper estimates for the logarithmic term.
\beqy
J'_{n,2i+1,2j+1}&=&\int_\fr{\O}{\d_{n,2j+1}}2\b_{2j+1}^2\fr{|y|^{\b_{2j+1}-2}}{(1+|y|^{\b_{2j+1}})^2}\wt\phi_{n,2j+1}(y)\mrm Pw_{n,2i+1}(\d_{n,2j+1}y)\mrm dy\\
&=&(-2\b_{2i+1}\log\d_{n,2i+1}+4\pi\b_{2i+1}H(0,0))\int_\fr{\O}{\d_{n,2j+1}}2\b_{2j+1}^2\fr{|y|^{\b_{2j+1}-2}}{(1+|y|^{\b_{2j+1}})^2}\wt\phi_{n,2j+1}(y)\mrm dy\\
&-&2\int_\fr{\O}{\d_{n,2j+1}}2\b_{2j+1}^2\fr{|y|^{\b_{2j+1}-2}}{(1+|y|^{\b_{2j+1}})^2}\wt\phi_{n,2j+1}(y)\log\l(1+\l(\fr{\d_{n,2j+1}}{\d_{n,2i+1}}\r)^{\b_{2i+1}}|y|^{\b_{2i+1}}\r)\mrm dy\\
&+&O\l(\d_{n,2j+1}\int_\fr{\O}{\d_{n,2j+1}}2\b_{2j+1}^2\fr{|y|^{\b_{2j+1}-1}}{(1+|y|^{\b_{2j+1}})^2}\l|\wt\phi_{n,2j+1}(y)\r|\mrm dy\r)\\
&+&O\l(\d_{n,2i+1}^{\b_{2i+1}}\int_\fr{\O}{\d_{n,2j+1}}2\b_{2j+1}^2\fr{|y|^{\b_{2j+1}-2}}{(1+|y|^{\b_{2j+1}})^2}\l|\wt\phi_{n,2j+1}(y)\r|\mrm dy\r)\\
&=&O\l(\log\fr{1}{|\la_n|}\r)\l|\int_\fr{\O}{\d_{n,2j+1}}2\b_{2j+1}^2\fr{|y|^{\b_{2j+1}-2}}{(1+|y|^{\b_{2j+1}})^2}\wt\phi_{n,2j+1}(y)\mrm dy\r|\\
&+&O\l(\l(\int_\fr{\O}{\d_{n,2i+1}}\fr{|y|^{(\b_{2j+1}-2)q}}{(1+|y|^{\b_{2j+1}})^{2q}}\log\l(1+\l(\fr{\d_{n,2j+1}}{\d_{n,2i+1}}\r)^{\b_{2i+1}}|y|^{\b_{2i+1}}\r)^q\r)^\fr{1}q\l\|\wt\phi_{n,2j+1}\r\|_\fr{q}{q-1}\r)\\
&+&O\l(\d_{n,2i+1}\l(\int_\fr{\O}{\d_{n,2i+1}}\fr{|y|^{(\b_{2j+1}-1)q}}{(1+|y|^{\b_{2j+1}})^{2q}}\r)^\fr{1}q\l\|\wt\phi_{n,2j+1}\r\|_\fr{q}{q-1}\r)+O\l(\d_{n,2i+1}^{\b_{2i+1}}\r)\\
&=&O(|\s_{n,2j+1}|)+O\l(\l(\fr{\d_{n,2j+1}}{\d_{n,2i+1}}\r)^{\b_{2i+1}}\l(\int_{B_\fr{\d_{n,2i+1}}{\d_{n,2j+1}}(0)}\fr{|y|^{(\b_{2j+1}+\b_{2i+1}-2)q}}{(1+|y|^{\b_{2j+1}})^{2q}}\mrm dy\r)^\fr{1}q\d_{n,2j+1}^{-2\l(1-\fr{1}q\r)}\|\phi_n\|\r)\\
&+&O\l(\l(\int_{\R^2\sm B_\fr{\d_{n,2i+1}}{\d_{n,2j+1}}(0)}\fr{|y|^{(\b_{2j+1}-2)q}}{(1+|y|^{\b_{2j+1}})^{2q}}\log\l(1+|y|^{\b_{2i+1}}\r)^q\mrm dy\r)^\fr{1}q\d_{n,2j+1}^{-2\l(1-\fr{1}q\r)}\|\phi_n\|\r)\\
&+&O\l(\d_{n,2i+1}\d_{n,2i+1}^{\min\l\{0,\b_{2i+1}+1-\fr{2}q\r\}}\d_{n,2j+1}^{-2\l(1-\fr{1}q\r)}\|\phi_n\|\r)+o(1)\\
&=&O\l(\l(\fr{\d_{n,2j+1}}{\d_{n,2i+1}}\r)^{\b_{2i+1}}\l(\fr{\d_{n,2j+1}}{\d_{n,2i+1}}\r)^{\min\l\{0,\b_{2i+1}-\b_{2j+1}+2\l(1-\fr{1}q\r)\r\}}\d_{n,2j+1}^{-2\l(1-\fr{1}q\r)}\r)\\
&+&O\l(\l(\fr{\d_{n,2j+1}}{\d_{n,2i+1}}\r)^{\fr{\b_{2j+1}}2+2\l(1-\fr{1}q\r)}\d_{n,2j+1}^{-2\l(1-\fr{1}q\r)}\r)+o(1)\\
&=&o(1).
\eeqy
In the other cases, some terms will vanish by the same arguments as before, but some others will not. To estimate the latter terms, we will use the convergence of $\wt\phi_{n,\ell}$ in $L_{\b_{2i+1}}$ and the following equalities, which can be proved by direct computation:
\beqy
&&\int_{\R^2}2\b_\ell^2\fr{|y|^{\b_\ell-2}}{(1+|y|^{\b_\ell}))^2}\fr{1-|y|^{\b_\ell}}{1+|y|^{\b_\ell}}\log\l(1+|y|^{\b_\ell}\r)\mrm dy=-2\pi\b_\ell;\\
&&\int_{\R^2}2\b_\ell^2\fr{|y|^{\b_\ell-2}}{(1+|y|^{\b_\ell}))^2}\fr{1-|y|^{\b_\ell}}{1+|y|^{\b_\ell}}\log|y|\mrm dy=-4\pi.
\eeqy
When $j=i$ we have
\beqy
J'_{n,2i+1,2i+1}&=&\int_\fr{\O}{\d_{n,2i+1}}2\b_{2i+1}^2\fr{|y|^{\b_{2i+1}-2}}{(1+|y|^{\b_{2i+1}})^2}\wt\phi_{n,2i+1}(y)\mrm Pw_{n,2i+1}(\d_{n,2i+1}y)\mrm dy\\
&=&(-2\b_{2i+1}\log\d_{n,2i+1}+4\pi\b_{2i+1}H(0,0))\int_\fr{\O}{\d_{n,2i+1}}2\b_{2i+1}^2\fr{|y|^{\b_{2i+1}-2}}{(1+|y|^{\b_{2i+1}})^2}\wt\phi_{n,2i+1}(y)\mrm dy\\
&+&\int_\fr{\O}{\d_{n,2i+1}}2\b_{2i+1}^2\fr{|y|^{\b_{2i+1}-2}}{(1+|y|^{\b_{2i+1}})^2}\wt\phi_{n,2i+1}(y)\log\l(1+|y|^{\b_{2i+1}}\r)\mrm dy\\
&+&O\l(\d_{n,2i+1}\int_\fr{\O}{\d_{n,2i+1}}2\b_{2i+1}^2\fr{|y|^{\b_{2i+1}-1}}{(1+|y|^{\b_{2i+1}})^2}\l|\wt\phi_{n,2i+1}(y)\r|\mrm dy\r)\\
&+&O\l(\d_{n,2i+1}^{\b_{2i+1}}\int_\fr{\O}{\d_{n,2i+1}}2\b_{2i+1}^2\fr{|y|^{\b_{2i+1}-2}}{(1+|y|^{\b_{2i+1}})^2}\l|\wt\phi_{n,2i+1}(y)\r|\mrm dy\r)\\
&=&-2\int_\fr{\O}{\d_{n,2i+1}}2\b_{2i+1}^2\fr{|y|^{\b_{2i+1}-2}}{(1+|y|^{\b_{2i+1}})^2}\wt\phi_{n,2i+1}(y)\log\l(1+|y|^{\b_{2i+1}}\r)\mrm dy+o(1)\\
&=&-2\mu_{2i+1}\int_{\R^2}2\b_{2i+1}^2\fr{|y|^{\b_{2i+1}-2}}{(1+|y|^{\b_{2i+1}})^2}\fr{1-|y|^{\b_{2i+1}}}{1+|y|^{\b_{2i+1}}}\log\l(1+|y|^{\b_{2i+1}}\r)+o(1)\\
&=&4\b_{2i+1}\mu_{2i+1}+o(1).
\eeqy
Similarly, if $j>i$,
\beqy
J'_{n,2i+1,2j+1}&=&\int_\fr{\O}{\d_{n,2j+1}}2\b_{2j+1}^2\fr{|y|^{\b_{2j+1}-2}}{(1+|y|^{\b_{2j+1}})^2}\wt\phi_{n,2j+1}(y)\mrm Pw_{n,2i+1}(\d_{n,2j+1}y)\mrm dy\\
&=&\l(-2\b_{2i+1}\log\d_{n,2j+1}+4\pi\b_{2i+1}H(0,0)\r)\int_\fr{\O}{\d_{n,2j+1}}2\b_{2j+1}^2\fr{|y|^{\b_{2j+1}-2}}{(1+|y|^{\b_{2j+1}})^2}\wt\phi_{n,2j+1}(y)\mrm dy\\
&-&2\b_{2i+1}\int_\fr{\O}{\d_{n,2j+1}}2\b_{2j+1}^2\fr{|y|^{\b_{2j+1}-2}}{(1+|y|^{\b_{2j+1}})^2}\wt\phi_{n,2j+1}(y)\log|y|\mrm dy\\
&-&2\int_\fr{\O}{\d_{n,2j+1}}2\b_{2j+1}^2\fr{|y|^{\b_{2j+1}-2}}{(1+|y|^{\b_{2j+1}})^2}\wt\phi_{n,2j+1}(y)\log\l(\l(\fr{\d_{n,2i+1}}{\d_{n,2j+1}}\r)^{\b_{2i+1}}\fr{1}{|y|^{\b_{2i+1}}}+1\r)\mrm dy\\
&+&O\l(\d_{n,2j+1}\int_\fr{\O}{\d_{n,2j+1}}2\b_{2j+1}^2\fr{|y|^{\b_{2j+1}-1}}{(1+|y|^{\b_{2j+1}})^2}\l|\wt\phi_{n,2j+1}(y)\r|\mrm dy\r)\\
&+&O\l(\d_{n,2i+1}^{\b_{2i+1}}\int_\fr{\O}{\d_{n,2j+1}}2\b_{2j+1}^2\fr{|y|^{\b_{2j+1}-2}}{(1+|y|^{\b_{2j+1}})^2}\l|\wt\phi_{n,2j+1}(y)\r|\mrm dy\r)\\
&=&-2\b_{2j+1}\mu_{2j+1}\int_{\R^2}2\b_{2j+1}^2\fr{|y|^{\b_{2j+1}-2}}{(1+|y|^{\b_{2j+1}})^2}\fr{1-|y|^{\b_{2j+1}}}{1+|y|^{\b_{2j+1}}}\log|y|\\
&+&O\l(\l(\int_\fr{\O}{\d_{n,2i+1}}\fr{|y|^{(\b_{2j+1}-2)q}}{(1+|y|^{\b_{2j+1}})^{2q}}\log\l(\l(\fr{\d_{n,2i+1}}{\d_{n,2j+1}}\r)^{\b_{2i+1}}\fr{1}{|y|^{\b_{2i+1}}}+1\r)^q\r)^\fr{1}q\l\|\wt\phi_{n,2j+1}\r\|_\fr{q}{q-1}\r)+o(1)\\
&=&8\pi\b_{2j+1}\mu_{2j+1}+O\l(\l(\int_{B_\fr{\d_{n,2i+1}}{\d_{n,2j+1}}(0)}\fr{|y|^{(\b_{2j+1}-2)q}}{(1+|y|^{\b_{2j+1}})^{2q}}\log\l(1+\fr{1}{|y|^{\b_{2i+1}}}\r)^q\mrm dy\r)^\fr{1}q\d_{n,2j+1}^{-2\l(1-\fr{1}q\r)}\|\phi_n\|\r)\\
&+&O\l(\l(\fr{\d_{n,2i+1}}{\d_{n,2j+1}}\r)^{\b_{2i+1}}\l(\int_{\R^2\sm B_\fr{\d_{n,2i+1}}{\d_{n,2j+1}}(0)}\fr{|y|^{(\b_{2j+1}-\b_{2i+1}-2)q}}{(1+|y|^{\b_{2j+1}})^{2q}}\mrm dy\r)^\fr{1}q\d_{n,2j+1}^{-2\l(1-\fr{1}q\r)}\|\phi_n\|\r)+o(1)\\
&=&8\pi\b_{2j+1}\mu_{2j+1}+O\l(\l(\fr{\d_{n,2i+1}}{\d_{n,2j+1}}\r)^{\fr{\b_{2j+1}}2-2\l(1-\fr{1}q\r)}\d_{n,2j+1}^{-2\l(1-\fr{1}q\r)}\r)\\
&+&O\l(\l(\fr{\d_{n,2i+1}}{\d_{n,2j+1}}\r)^{\b_{2i+1}}\l(\fr{\d_{n,2i+1}}{\d_{n,2j+1}}\r)^{\min\l\{0,\b_{2j+1}-\b_{2i+1}-2\l(1-\fr{1}q\r)\r\}}\d_{n,2j+1}^{-2\l(1-\fr{1}q\r)}\r)+o(1)\\
&=&8\pi\b_{2j+1}\mu_{2j+1}+o(1)
\eeqy
$J''_{n,2i+1}$ vanishes because, by Lemma \ref{pw}, $\|\mrm Pw_{n,\ell}\|_\infty=O\l(\log\fr{1}{|\la_n|}\r)$, therefore
$$\l|J''_{n,2i+1}\r|\le\|\psi_{n,1}\|_1\|\mrm Pw_{n,2i+1}\|_\infty\le C\log\fr{1}{|\la_n|}\|\psi_n\|_p\us{n\to+\infty}\to0.$$
Finally, Step $3$ gives
$$J'''_{n,2i+1}=\int_\fr{\O}{\d_{n,2i+1}}2\b_{2i+1}^2\fr{|y|^{\b_{2i+1}-2}}{(1+|y|^{\b_{2i+1}})^2}\wt\phi_{n,2i+1}(y)\mrm dy=\fr{\s_{n,2i+1}}{\log\fr{1}{|\la_n|}}\us{n\to+\infty}\to0$$
Putting all these estimates together, repeating the computations for even indexes and passing to the limit gives
\beqy
4\pi\b_{2i+1}\mu_{2i+1}+8\pi\sum_{j=i+1}^{\l[\fr{k-1}2\r]}\b_{2j+1}\mu_{2j+1}-4\pi a\sum_{j=i}^{\l[\fr{k-2}2\r]}\b_{2j+2}\mu_{2j+2}&=&0\\
4\pi\b_{2i+2}\mu_{2i+2}+8\pi\sum_{j=i+1}^{\l[\fr{k-2}2\r]}\b_{2j+2}\mu_{2j+2}-4\pi b\sum_{j=i+1}^{\l[\fr{k-1}2\r]}\b_{2j+1}\mu_{2j+1}&=&0,
\eeqy
from which we get $\mu_\ell=0$ for all $j$'s.\\

\item[\emph{Step $5$:}]\emph{$\phi_n\us{n\to+\infty}\to0$ in $L^\infty(\O)^2$.}\\
We fix $x\in\O$ and we estimate $\phi_{n,i}(x)$, using Green's representation formula. We provide the estimate only for $i=1$:
\beqy
|\phi_{n,1}(x)|&=&\l|\sum_{j=0}^{\l[\fr{k-1}2\r]}\int_\O G(x,y)2\b_{2j+1}^2\fr{\d_{n,2j+1}^{\b_{2j+1}}|y|^{\b_{2j+1}-2}}{\l(\d_{n,2j+1}^{\b_{2j+1}}+|y|^{\b_{2j+1}}\r)^2}\phi_{n,1}(y)\mrm dy\r.\\
&-&\l.\fr{a}2\sum_{j=0}^{\l[\fr{k-2}2\r]}\int_\O G(x,y)2\b_{2j+2}^2\fr{\d_{n,2j+2}^{\b_{2j+2}}|y|^{\b_{2j+2}-2}}{\l(\d_{n,2j+2}^{\b_{2j+2}}+|y|^{\b_{2j+2}}\r)^2}\phi_{n,2}(y)\mrm dy+\int_\O G(x,y)\psi_{n,1}(y)\mrm dy\r|\\
&\le&\sum_{i=1}^2\sum_{j=0}^{\l[\fr{k-i}2\r]}\l|\int_\O G(x,y)\fr{\d_{n,2j+i}^{\b_{2j+i}}|y|^{\b_{2j+i}-2}}{\l(\d_{n,2j+i}^{\b_{2j+i}}+|y|^{\b_{2j+i}}\r)^2}\phi_{n,i}(y)\mrm dy\r|+\l|\int_\O G(x,y)\psi_{n,1}(y)\mrm dy\r|\\
&\le&\sum_{j=0}^k\l|\int_\fr{\O}{\d_{n,\ell}}G(x,\d_{n,\ell}z)\fr{|z|^{\b_\ell-2}}{\l(1+|z|^{\b_\ell}\r)^2}\wt\phi_{n,\ell}(z)\mrm dz\r|+\sup_{x\in\O}\|G(x,\cd)\|_\fr{p}{p-1}\|\psi_n\|_p\\
&\le&\sum_{j=0}^k\l|\ub{\int_\fr{\O}{\d_{n,\ell}}\log|x-\d_{n,\ell}z|\fr{|z|^{\b_\ell-2}}{\l(1+|z|^{\b_\ell}\r)^2}\wt\phi_{n,\ell}(z)\mrm dz}_{:=K'_{n,\ell}}\r|+\sum_{j=0}^k\ub{\int_\fr{\O}{\d_{n,\ell}}|H(x,\d_{n,\ell}z)|\fr{|z|^{\b_\ell-2}}{\l(1+|z|^{\b_\ell}\r)^2}\l|\wt\phi_{n,\ell}(z)\r|\mrm dz}_{:=K''_{n,\ell}}+o(1).
\eeqy
To estimate $K''_{n,\ell}$ we apply some weighted Sobolev inequalities to $\wt\phi_{n,\ell}$: since it is bounded in $H_{\b_\ell}\l(\R^2\r)$ and tends to $0$ in $L_{\b_\ell}\l(\R^2\r)$, then for any $q\ge2$
$$\int_\fr{\O}{\d_{n,\ell}}\fr{|z|^{\b_\ell-2}}{(1+|z|^{\b_\ell})^2}\l|\wt\phi_{n,\ell}(z)\r|^q\mrm dz\us{n\to+\infty}\to0.$$
Therefore, for a suitable $q$,
\beqy
K''_{n,\ell}&\le&(H(0,0)+|x|)\int_\fr{\O}{\d_{n,\ell}}\fr{|z|^{\b_\ell-2}}{\l(1+|z|^{\b_\ell}\r)^2}\l|\wt\phi_{n,\ell}(z)\r|\mrm dz+\d_{n,\ell}\int_\fr{\O}{\d_{n,\ell}}\fr{|z|^{\b_\ell-1}}{\l(1+|z|^{\b_\ell}\r)^2}\l|\wt\phi_{n,\ell}(z)\r|\mrm dz\\
&\le&C\l(\int_\fr{\O}{\d_{n,\ell}}\fr{|z|^{\b_\ell-2}}{\l(1+|z|^{\b_\ell}\r)^2}\mrm dz\r)^\fr{1}2\l(\int_\fr{\O}{\d_{n,\ell}}\fr{|z|^{\b_\ell-2}}{\l(1+|z|^{\b_\ell}\r)^2}\l|\wt\phi_{n,\ell}(z)\r|^2\mrm dz\r)^\fr{1}2\\
&+&\d_{n,\ell}\l(\int_\fr{\O}{\d_{n,\ell}}\fr{|z|^{\b_\ell-2+\fr{q}{q-1}}}{\l(1+|z|^{\b_\ell}\r)^2}\mrm dz\r)^{1-\fr{1}q}\l(\int_\fr{\O}{\d_{n,\ell}}\fr{|z|^{\b_\ell-2}}{\l(1+|z|^{\b_\ell}\r)^2}\l|\wt\phi_{n,\ell}(z)\r|^q\mrm dz\r)^\fr{1}q\\
&\le&C\l\|\wt\phi_{n,\ell}\r\|_{L_{\b_\ell}\l(\R^2\r)}+\d_{n,\ell}^{\min\l\{1,\b_\ell\l(1-\fr{1}q\r)\r\}}o(1)\\
&\us{n\to+\infty}\to&0.
\eeqy
To deal with $K'_{n,\ell}$ we use that $\s_{n,\ell}\us{n\to+\infty}\to0$ (see Step $3$) and that $\d_{n,\ell}$ is given by powers of $\la_{n,i}$; in particular, we will distinguish whether $|x|$ is smaller or larger than $\d_{n,\ell}$:
\beqy
|K'_{n,\ell}|&\le&\l|\int_\fr{\O}{\d_{n,\ell}}\log\l|\fr{x-\d_{n,\ell}z}{\max\{\d_{n,\ell},|x|\}}\r|\fr{|z|^{\b_\ell-2}}{\l(1+|z|^{\b_\ell}\r)^2}\wt\phi_{n,\ell}(z)\mrm dz\r|+|\log\max\{\d_{n,\ell},|x|\}|\l|\int_\fr{\O}{\d_{n,\ell}}\fr{|z|^{\b_\ell-2}}{\l(1+|z|^{\b_\ell}\r)^2}\wt\phi_{n,\ell}(z)\mrm dz\r|\\
&\le&\l(\int_\fr{\O}{\d_{n,\ell}}\l|\log\fr{\l|\fr{x}{\d_{n,\ell}}-z\r|}{\max\l\{1,\l|\fr{x}{\d_{n,\ell}}\r|\r\}}\r|^2\fr{|z|^{\b_\ell-2}}{\l(1+|z|^{\b_\ell}\r)^2}\mrm dz\r)^\fr{1}2\l(\int_\fr{\O}{\d_{n,\ell}}\fr{|z|^{\b_\ell-2}}{\l(1+|z|^{\b_\ell}\r)^2}\l|\wt\phi_{n,\ell}(z)\r|^2\mrm dz\r)^\fr{1}2\\
&+&\fr{\max\l\{\log\fr{1}{\d_{n,\ell}},\log\mrm{diam}\O\r\}}{\log\fr{1}{|\la_n|}}|\s_{n,\ell}|\\
&\le&\l(\ub{\int_{\R^2}\l|\log\fr{|z'|}{\max\l\{1,\l|\fr{x}{\d_{n,\ell}}\r|\r\}}\r|^2\fr{\l|z'-\fr{x}{\d_{n,\ell}}\r|^{\b_\ell-2}}{\l(1+\l|z'-\fr{x}{\d_{n,\ell}}\r|^{\b_\ell}\r)^2}\mrm dz'}_{K'''_\ell\l(\fr{x}{\d_{n,\ell}}\r)}\r)^\fr{1}2o(1)+o(1);
\eeqy
The claim will follow by showing that $K'''_\ell(x')$ is uniformly bounded for $x'\in\R^2$.\\
Taking a cue from \cite{ci} (Lemma $1.1$), we split the integral in the ball of radius $2\max\{1,|x'|\}$ and its complementary: in the ball, we just apply a H\"older inequality with suitable exponents and then a dilatation; in its exterior, we use the monotonicity of the logarithm and the fact that $x'$ it is somehow negligible with respect to $z'$:
$$\fr{|z'|}2\le|z'|-\max\{1,|x'|\}\le|z'|-|x'|\le|z'-x'|\le|z'|+|x'|\le|z'|+\max\{1,|x'|\}\le\fr{3}2|z'|.$$
We get, for a suitable $q>1$:
\beqy
K'''_\ell(x)&=&\int_{B_{2\max\{1,|x'|\}}(0)}\l|\log\fr{|z'|}{\max\{1,|x'|\}}\r|^2\fr{|z'-x'|^{\b_\ell-2}}{(1+|z'-x'|^{\b_\ell})^2}\mrm dz'\\
&+&\int_{\R^2\sm B_{2\max\{1,|x'|\}}(0)}\l|\log\fr{|z'|}{\max\{1,|x'|\}}\r|^2\fr{|z'-x'|^{\b_\ell-2}}{(1+|z'-x'|^{\b_\ell})^2}\mrm dz'\\
&\le&\l(\int_{B_{2\max\{1,|x'|\}}(0)}\l|\log\fr{|z'|}{\max\{1,|x'|\}}\r|^\fr{2q}{q-1}\mrm dz'\r)^{1-\fr{1}q}\l(\int_{\R^2}\fr{|z'-x'|^{(\b_\ell-2)q}}{(1+|z'-x'|^{\b_\ell})^{2q}}\mrm dz'\r)^\fr{1}q\\
&+&C\int_{\R^2\sm B_{2\max\{1,|x'|\}}(0)}\l|\log\fr{|z'|}{\max\{1,|x'|\}}\r|^2\fr{|z'|^{\b_\ell-2}}{\l(1+|\fr{z'}2|^{\b_\ell}\r)^2}\mrm dz'\\
&\le&C\l(\int_{B_2(0)}|\log|y'||^\fr{2q}{q-1}\mrm dy'\r)^{1-\fr{1}q}+C\int_{\R^2\sm B_2(0)}|\log|z'||^2\fr{|z'|^{\b_\ell-2}}{\l(1+|\fr{z'}2|^{\b_\ell}\r)^2}\mrm dz'\\
&\le&C.
\eeqy\

\item[\emph{Step $6$:}]\emph{A contradiction arises.}\\
We multiply each equation of \eqref{phn} by the respective $\phi_{n,i}$ and we sum the two of them. We get:
\beqy
1&=&\int_\O|\n\phi_{n,1}|^2+\int_\O|\n\phi_{n,2}|^2\\
&=&\sum_{j=0}^{\l[\fr{k-1}2\r]}\int_\O2\b_{2j+1}^2\fr{\d_{n,2j+1}^{\b_{2j+1}}|\cd|^{\b_{2j+1}-2}}{\l(\d_{n,2j+1}^{\b_{2j+1}}+|\cd|^{\b_{2j+1}}\r)^2}\phi_{n,1}^2-\fr{a}2\sum_{j=0}^{\l[\fr{k-2}2\r]}\int_\O2\b_{2j+2}^2\fr{\d_{n,2j+2}^{\b_{2j+2}}|\cd|^{\b_{2j+2}-2}}{\l(\d_{n,2j+2}^{\b_{2j+2}}+|\cd|^{\b_{2j+2}}\r)^2}\phi_{n,1}\phi_{n,2}+\int_\O\psi_{n,1}\phi_{n,1}\\
&+&\sum_{j=0}^{\l[\fr{k-2}2\r]}\int_\O2\b_{2j+2}^2\fr{\d_{n,2j+2}^{\b_{2j+2}}|\cd|^{\b_{2j+2}-2}}{\l(\d_{n,2j+2}^{\b_{2j+2}}+|\cd|^{\b_{2j+2}}\r)^2}\phi_{n,2}^2-\fr{b}2\sum_{j=0}^{\l[\fr{k-1}2\r]}\int_\O2\b_{2j+1}^2\fr{\d_{n,2j+1}^{\b_{2j+1}}|\cd|^{\b_{2j+1}-2}}{\l(\d_{n,2j+1}^{\b_{2j+1}}+|\cd|^{\b_{2j+1}}\r)^2}\phi_{n,1}\phi_{n,2}+\int_\O\psi_{n,2}\phi_{n,2}\\
&\le&C\sum_{j=0}^k\int_\O\fr{\d_{n,\ell}^{\b_\ell}|\cd|^{\b_\ell-2}}{\l(\d_{n,\ell}^{\b_\ell}+|\cd|^{\b_\ell}\r)^2}\|\phi_n\|^2_{L^\infty(\O)}+\|\psi_n\|_p\|\phi_n\|_\fr{p}{p-1}\\
&\le&C\l(\|\phi_n\|^2_{L^\infty(\O)}+\|\psi_n\|_p\|\phi_n\|_\fr{p}{p-1}\r)\\
&\us{n\to+\infty}\to&0;
\eeqy
which is a contradiction.
\eite

\epf\

\appendix\

\section{Appendix}\label{appe}\

We prove here a classification result for entire solutions of a scalar linearized problem.\\

\bprop\label{intere}${}$\\
Assume $\a>0$, $m\in\N$ and $\fr{\a}{ m}\nin2\N$. Then, any solution $\phi$ of
\bequ\label{solint}
\l\{\bll-\D\phi=2\a^2\fr{|\cd|^{\a-2}}{\l(1+|\cd|^\a\r)^2}\phi&\tx{in }\R^2\\\int_{\R^2}|\n\phi|^2<+\infty\\\phi\l(e^{\fr{2\pi}m\iota}\cd\r)=\phi\earr\r.
\eequ
satisfies, for some $\mu\in\R$,
$$\phi=\mu\fr{1-|\cd|^\a}{1+|\cd|^\a}.$$
\eprop\

\bpf${}$\\
We argue as Baraket and Pacard do in \cite{barpac}, Proposition $1$ where the case $\a=2$ is covered (see also Del Pino, Esposito and Musso \cite{dem2}).\\
By writing any solution $\phi$ of \eqref{solint} as a Fourier decomposition
$$\phi(x)=\sum_{n\in\Z}\phi_n(|x|)e^{in\t},$$
we see that each of the $\phi_n$ solves the following o.d.e.
\bequ\label{phin}
\pa_\rho^2\phi_n(\rho)+\fr{1}\rho\pa_\rho\phi_n(\rho)-\fr{n^2}{\rho^2}\phi_n(\rho)+\fr{2\a^2\rho^{\a-2}}{\l(1+\rho^\a\r)^2}\phi_n(\rho).
\eequ
Integrating by parts, $\phi_n$ must satisfy
$$\int_0^{+\infty}\l(\l|\pa_\rho\phi_n(\rho)\r|^2+\l(\fr{n^2}{\rho^2}-\fr{2\a^2\rho^{\a-2}}{\l(1+\rho^\a\r)^2}\r)\phi_n(\rho)^2\r)\rho\mrm d\rho=0;$$
since $\fr{n^2}{\rho^2}-\fr{2\a^2\rho^{\a-2}}{\l(1+\rho^\a\r)^2}\ge\fr{1}{\rho^2}\l(n^2-\fr{\a^2}2\r)$, we must have $\phi_n\eq0$ for $|n|\ge\fr{\a}{\sqrt2}$. In particular, $\phi$ is a finite combination of the $\phi_n$'s.\\
It is easy to see that each solution of \eqref{phin} is a linear combination of the fundamental solutions
$$\phi_{n,+}(\rho)=\rho^n\fr{\a+2n-(\a-2n)\rho^\a}{1+\rho^\a}\q\q\q\phi_{n,-}(\rho)=\rho^{-n}\fr{\a-2n-(\a+2n)\rho^\a}{1+\rho^\a}.$$
Since we are looking for bounded solutions of \eqref{solint}, here we are allowed to take only bounded solution of \eqref{phin}.\\
If $\a$ is not an even integer, the condition is satisfied only by $\phi_{0,+}(\rho)=\phi_{0,-}(\rho)=\fr{1-\rho^\a}{1+\rho^\a}$, hence $\phi(x)=\phi_0(|x|)$ is an integer multiple of its and the Proposition is proved.\\
On the other hand, if $\a\in2\N$, then $\phi_{\fr{\a}2,+}(\rho)=2\a\fr{\rho^\fr{\a}2}{1+\rho^\a}$ is also allowed, therefore in this case $\phi(x)$ is a combination of the following functions:
$$\phi_0(|x|)=\fr{1-|x|^\a}{1+|x|^\a}\q\q\q\fr{1}{2\a}\phi_{\fr{\a}2}(|x|)\cos\l(\fr{\a}2\t\r)=\fr{|x|^\fr{\a}2}{1+|x|^\a}\cos\l(\fr{\a}2\t\r)\q\q\q\fr{1}{2\a}\phi_{\fr{\a}2}(|x|)\sin\l(\fr{\a}2\t\r)=\fr{|x|^\fr{\a}2}{1+|x|^\a}\sin\l(\fr{\a}2\t\r).$$
Anyway, the latter two functions do not satisfy the symmetry requirement if $m$ is as in the assumptions, therefore $\phi$ must again be a multiple of $\phi_0(|x|)$. The proof is completed.
\epf\

\bibliography{2x2}
\bibliographystyle{abbrv}

\end{document}